# A FORWARD–BACKWARD STOCHASTIC ALGORITHM FOR QUASI-LINEAR PDES[1]

By François Delarue and Stéphane Menozzi

*Université Paris VII*

We propose a time-space discretization scheme for quasi-linear parabolic PDEs. The algorithm relies on the theory of fully coupled forward–backward SDEs, which provides an efficient probabilistic representation of this type of equation. The derivated algorithm holds for strong solutions defined on any interval of arbitrary length. As a bypass product, we obtain a discretization procedure for the underlying FBSDE. In particular, our work provides an alternative to the method described in [Douglas, Ma and Protter (1996) *Ann. Appl. Probab.* **6** 940–968] and weakens the regularity assumptions required in this reference.

**1. Introduction.** Introduced first by Antonelli [1] and then by Ma, Protter and Yong [14], forward–backward stochastic differential equations (FBSDEs in short) provide an extension of the Feynman–Kac representation to a certain class of quasi-linear parabolic PDEs. These equations also appear in a large number of application fields such as the Hamiltonian formulation of control problems or the option hedging problem with large investors in financial mathematics (i.e., when the wealth or strategy of an agent has an impact on the volatility). We refer to the monograph of Ma and Yong [15] for details and further applications.

1.1. *FBSDE theory and discretization algorithm.*

*Connection between FBSDEs and quasi-linear parabolic PDEs.* Consider a probability space $(\Omega, \mathcal{F}, \mathbb{P})$ endowed with a $d$-dimensional Brownian motion $(B_t)_{t \in [0,T]}$, where $T$ denotes an arbitrarily prescribed positive real. For a given initial condition $x_0 \in \mathbb{R}^d$, a forward–backward SDE strongly couples

Received September 2004; revised June 2005.

[1]Supported in part by the CMAP, Ecole Polytechnique.

a diffusion process $U$ to the solution $(V, W)$ of a backward SDE (as defined in the earlier work of Pardoux and Peng [20]):

$$\forall t \in [0, T] \qquad U_t = x_0 + \int_0^t b(U_s, V_s, W_s)\, ds + \int_0^t \sigma(U_s, V_s)\, dB_s,$$

(E)

$$V_t = H(U_T) + \int_t^T f(U_s, V_s, W_s)\, ds - \int_t^T W_s\, dB_s.$$

In this paper the coefficients $b$, $f$, $\sigma$ and $H$ are deterministic (and, for simplicity, also time independent). In this case, Ma, Protter and Yong [14], Pardoux and Tang [21] and Delarue [6] have investigated in detail the link with the following quasi-linear PDE on $[0, T[ \times \mathbb{R}^d$:

$$\partial_t u(t, x) + \langle b(x, u(t, x), \nabla_x u(t, x)\sigma(x, u(t, x))), \nabla_x u(t, x)\rangle$$

(ℰ)

$$+ \tfrac{1}{2}\mathrm{tr}(a(x, u(t, x))\nabla_{x,x}^2 u(t, x))$$

$$+ f(x, u(t, x), \nabla_x u(t, x)\sigma(x, u(t, x))) = 0,$$

$$u(T, x) = H(x),$$

with $a(x, y) = (\sigma\sigma^*)(x, y)$, $(x, y) \in \mathbb{R}^d \times \mathbb{R}$.

*A probabilistic numerical method for FBSDEs and quasi-linear PDEs.* This paper aims to derive from the probabilistic theory of FBSDEs a completely tractable algorithm to approximate the solution of equation (ℰ). As a bypass product, the procedure also provides a discretization of the triple $(U, V, W)$.

Most of the available numerical methods proposed so far are purely analytic and involve finite-difference or finite-element techniques to approximate the solution $u$ of (ℰ). For example, the discretization procedure for FBSDEs of type (E), given in [10], consists in discretizing first the PDE (ℰ) and then in deriving an approximation of the underlying FBSDE.

At the opposite, we propose in this paper to derive from the FBSDE representation a numerical scheme for quasi-linear equations of type (ℰ). This strategy finds its origin in the earlier work of Chevance [5], who introduced a time-space discretization scheme in the decoupled or so-called "pure backward" case. In this latter frame, the coefficients $b$ and $\sigma$ do not depend on $V$ and $W$ and the forward equation reduces to a classical SDE. The process $U$ then appears as an "objective diffusion." Note in this particular case that the time-space discretization scheme and the specific form of the system (E) permit to use a standard "dynamic programming principle."

From a numerical point of view, two other kinds of approaches have been developed in the backward case. The first one is based on Monte Carlo simulations and Malliavin integration by parts; see [4]. The other one relies on quantization techniques for a discretization scheme of the underlying forward equation. Quantization consists in approximating a random variable by a suitable discrete law. It provides a cheap and numerically efficient alternative to usual Monte Carlo methods to estimate expectations. In the works of Bally and Pagès [2] or Bally, Pagès and Printems [3] on American options, the key idea is to perform an optimal quantization procedure of a discretized version of the underlying diffusion process in order to compute *once for all* by a Monte Carlo method the corresponding semi-group. Then, the second step consists in doing a dynamic programming descent. For other applications of quantization, we refer to the works of Pagès, Pham and Printems [18] or Pagès and Printems [19].

*Discretization strategy.* In the coupled case, or quasi-linear framework, the diffusion $U$ is not "objective" anymore. Indeed, due to the strong nonlinearity of the equation $(\mathcal{E})$, the coefficients of the underlying forward diffusion depend on the solution and on its gradient.

In particular, we cannot quantify a discretization scheme of the diffusion process as explained above. This is well understood: without approximating $u$, we do not have any a priori knowledge of the optimal shape of the associated grid. Hence, we just focus on the quantization of the Brownian increments appearing in the forward SDE and then choose to define the approximate diffusion on a sequence of truncated $d$-dimensional Cartesian grids. Note that the discretization procedure of $U$ is now coupled to the approximation procedure of $(u, \nabla_x u)$ [denoted in a generic way by $(\bar{u}, \bar{v})$] which is computed along the same sequence of grids. The time-space discretization scheme allows to define $(\bar{u}, \bar{v})$ and the approximations of the transitions of $U$ in order to recover a kind of "dynamic programming principle." Consider indeed a given regular time mesh $(t_i = ih)_{i \in \{0,\dots,N\}}$ of $[0,T]$, $h$ being the step size. To every discretization time $t_i$, associate a spatial Cartesian grid $\mathcal{C}_i \equiv \{(x_k^i)_{k \in \mathcal{I}_i}\}$, $\mathcal{I}_i \subset \mathbb{N}^*$, such that $\forall i \in \{0,\dots,N-1\}$, $\mathcal{C}_i \subset \mathcal{C}_{i+1}$. Starting from $t_N = T$ for which the solution of $(\mathcal{E})$ and its gradient are known, the transition of $U$ from $t_i$ to $t_{i+1}$, $i \in \{0,\dots,N-1\}$, is then updated iteratively through the Brownian quantized increments and through the values of $\bar{u}(t_{i+1}, \cdot)$ and $\bar{v}(t_{i+1}, \cdot)$ on the grid $\mathcal{C}_{i+1}$. This permits to express the approximation $\bar{u}(t_i, \cdot)$ through a discretized version of the Feynman–Kac formula.

At this stage, it remains to specify the way we update the approximation of the gradient of the solution $u$. We mention actually that the strategy aims to approximate the product $\nabla_x u(t_k, \cdot) \sigma(\cdot, u(t_k, \cdot))$ instead of $\nabla_x u(t_k, \cdot)$ itself. This explains the specific writing of the PDE $(\mathcal{E})$. We then proceed in two different steps. A first approximation is performed through a martingale increment procedure as done in the discretization scheme of BSDEs explained in [4], or as used in [3]. A second step consists in quantizing the Gaussian increments appearing in the former representation. This is an alternative solution to the usual techniques based on Monte Carlo simulations or on Malliavin integration by parts as employed in [4]. Of course, if the matrix $\sigma\sigma^*$ is nondegenerate, the strategy still applies, up to an inversion procedure, to coefficients of the form $(b, f)(x, u(t, x), \nabla_x u(t, x))$.

*Extra references.* Some of the preliminaries of our approach can be found in [17] in the specific case where $(b, f)(x, u(t, x), \nabla_x u(t, x)\sigma(t, x, u(t, x)))$ reduces to $(b, f)(x, u(t, x))$. Note, however, that the proof of the convergence of the underlying numerical scheme proposed in this reference just holds for so-called "equations with small parameter" (i.e., with a small diffusion matrix). Generally speaking, the authors have then to control the regularity properties of the solution of the transport problem associated to the equation ($\mathcal{E}$) [i.e., the same equation as ($\mathcal{E}$), but without any second-order terms]. Without discussing in detail the basic assumptions made in our paper, note that no condition of this type appears in the sequel: in particular, the matrix $a$ is assumed to be uniformly elliptic. Hence, we feel that the work of Milstein and Tretyakov [17] applies to a different framework than ours. For this reason, we avoid any further comparisons between both situations. Add finally, for the sake of completeness, that Makarov [16] has successfully applied the strategy of Milstein and Tretyakov [17] to the case $(b, f) \equiv (b, f)(x, u(t, x), \nabla_x u(t, x)\sigma(t, x, u(t, x)))$ under suitable smoothness properties on the coefficients. Of course, the small parameter condition is then still necessary.

## 1.2. *Novelties brought by the paper.*

*A purely probabilistic point of view.* The proof of the convergence of our algorithm is somehow the first to be essentially of probabilistic nature, since we are able to adapt the usual stability techniques of BSDE theory to the discretized framework. Note, in particular, that we follow the proof of uniqueness in the *four step scheme* given in [14] to handle the strong coupling between the forward and backward components.

In the discretized framework, the gradient terms appearing in $b$ and $f$ bring additional difficulties. Indeed, our gradient approximation does not appear as a representation process given by the martingale representation theorem as the process $W$ in ($E$). In particular, the strategy introduced by Pardoux and Peng [20] to estimate the $L^2$ norm of $W$ over $[0, T]$ fails in the discretized setting. We then propose a specific probabilistic strategy to overcome this deep trouble and thus to handle the nonlinearities of order one, see Sections 3.3 and 9.3 for details.

*Convergence under weak assumptions.* In [10], the authors handle the gradient terms by working under smoothness assumptions that allow them to study the gradient of $u$ as the solution of the differentiated PDE.

Our strategy permits to avoid to differentiate the PDE and thus to really weaken the assumptions required both on the coefficients of $(E)$ and on the smoothness of the solution $u$ of $(\mathcal{E})$ in the above reference. In the previous paper, the coefficients are assumed to be smoothly differentiable and bounded. We just suppose that they are Lipschitz continuous and bounded in $x$. In [10], the solution $u$ of $(\mathcal{E})$ is at least bounded in $C^{2+\alpha/2,4+\alpha}([0,T] \times \mathbb{R}^d), \alpha \in ]0,1[$. In our paper we only impose $u$ to belong to $C^{1,2}([0,T] \times \mathbb{R}^d)$ with bounded derivatives of order one in $t$ and one and two in $x$.

*A completely tractable algorithm.* Furthermore, in [10], the authors always take into consideration the case of infinite spatial grids. This turns out to be simpler for the convergence analysis, anyhow it does not provide in all generality a fully implementable algorithm. We discuss the impact of the truncation of the grids and analyze its contribution in the error.

Finally, a linear interpolation procedure is also used in [10] to define the algorithm. This can be heavy in large dimension. The algorithm we propose allows to define the approximate solution only at the nodes of the spatial grid. In this way, we feel that our method is simpler to implement and numerically cheaper. Note, moreover, that we avoid the inversion of large linear systems associated to "usual" numerical analysis techniques.

1.3. *Organization of the paper.* In Section 2 we detail general assumption and notation, as well as several smoothness properties of the solution $u$ of $(\mathcal{E})$. We also specify the connection between the FBSDE $(E)$ and the quasi-linear PDE $(\mathcal{E})$. Section 3 explains the main algorithmic choices. We present, in particular, the various steps that led us to the current discretization scheme. The main results are stated and discussed in Section 4. In particular, we give an estimate of the speed of convergence of the algorithm. As a probabilistic counterpart, we estimate the difference between the approximating processes and the initial solution $(U,V,W)$ of $(E)$. Numerical examples are presented in Section 5.

The end of the paper is then mainly devoted to the proof of the convergence results. The proof is divided into three parts. Various a priori controls of the discrete objects are stated and proved in Section 6. In Section 7 we adapt the FBSDE machinery to our setting to prove a suitable stability property. Section 8 is then devoted to the last step of the proof and, more precisely, to a specific refinement of Gronwall's lemma. In order to be concise, we sometimes only sketch the proofs. They are presented in detail in the electronic version Delarue and Menozzi [9].

As a conclusion, we compare in Section 9 our strategy to other methods and explain some technical points that motivated the choice of our current algorithm. We also indicate further conceivable extensions.

**2. Nonlinear Feynman–Kac formula.** In this section we first give the assumptions on the coefficients of the FBSDE and then briefly recall the connection with quasi-linear PDEs. As detailed later, under these assumptions, the underlying PDE admits a unique strong solution, whose partial derivatives of order one in $t$ and one and two in $x$ are controlled on the whole domain by known parameters. For the sake of simplicity, we also assume that the coefficients do not depend on time.

2.1. *Coefficients of the equation.* For a given $d \in \mathbb{N}^*$, we consider the coefficients $b : \mathbb{R}^d \times \mathbb{R} \times \mathbb{R}^d \to \mathbb{R}^d$, $f : \mathbb{R}^d \times \mathbb{R} \times \mathbb{R}^d \to \mathbb{R}$, $\sigma : \mathbb{R}^d \times \mathbb{R} \to \mathbb{R}^{d \times d}$, $H : \mathbb{R}^d \to \mathbb{R}$.

ASSUMPTION (A). We say that the functions $b$, $f$, $H$ and $\sigma$ satisfy Assumption (A) if they are bounded in space and have at most linear growth in the other variables, are uniformly Lipschitz continuous w.r.t. all the variables, $a = \sigma \sigma^*$ is uniformly elliptic and $H$ is bounded in $C^{2+\alpha}(\mathbb{R}^d)$.

*From now on, Assumption (A) is in force.*

2.2. *Forward–backward SDE.* Consider now a given $T > 0$ and a probability space $(\Omega, \mathcal{F}, \mathbb{P})$ endowed with a Brownian motion $(B_t)_{0 \leq t \leq T}$ whose natural filtration, augmented with $\mathbb{P}$ null sets, is denoted by $\{\mathcal{F}_t\}_{0 \leq t \leq T}$.

Fix an initial condition $x_0 \in \mathbb{R}^d$ and recall (see [14] and [6]) that there exists a unique progressively measurable triple $(U, V, W)$, with values in $\mathbb{R}^d \times \mathbb{R} \times \mathbb{R}^d$, such that $\mathbb{E} \sup_{t \in [0,T]}(|U_t|^2 + |V_t|^2) < +\infty$, $\mathbb{E} \int_0^T |W_t|^2 \, dt < +\infty$, and which satisfies $\mathbb{P}$ almost surely the couple of equations (E).

2.3. *Quasi-linear PDE.* Thanks to [13], Chapter VI, Theorem 4.1, and to [14] (up to a regularization procedure of the coefficients), we claim that ($\mathcal{E}$) admits a solution $u \in \mathcal{C}^{1,2}([0,T] \times \mathbb{R}^d, \mathbb{R})$ satisfying the following:

THEOREM 2.1. *There exists a constant $C_{2.1}$, only depending on $T$ and on known parameters deriving from Assumption (A), such that $\forall (t,x) \in [0,T] \times \mathbb{R}^d$,*

$$|u(t,x)| + |\nabla_x u(t,x)| + |\nabla_{x,x}^2 u(t,x)| + |\partial_t u(t,x)|$$
$$+ \sup_{t' \in [0,T], t \neq t'} [|t - t'|^{-1/2} |\nabla_x u(t,x) - \nabla_x u(t',x)|] \leq C_{2.1}.$$

*Moreover, $u$ is unique in the class of functions $\widetilde{u} \in \mathcal{C}([0,T] \times \mathbb{R}^d, \mathbb{R}) \cap \mathcal{C}^{1,2}([0,T[ \times \mathbb{R}^d, \mathbb{R})$ which satisfy $\sup_{(t,x) \in [0,T[ \times \mathbb{R}^d}(|\widetilde{u}(t,x)| + |\nabla_x \widetilde{u}(t,x)|) < +\infty$.*

From [6, 14, 21], the FBSDE $(E)$ is connected with the PDE $(\mathcal{E})$. Set $\forall (t,x) \in [0,T[\times \mathbb{R}^d$, $v(t,x) = \nabla_x u(t,x)\sigma(x,u(t,x))$. The relationship between $(E)$ and $(\mathcal{E})$ can be summed up as follows: $\forall t \in [0,T]$,

(2.1)
$$V_t = u(t,U_t), \qquad W_t = v(t,U_t),$$

$$V_t = \mathbb{E}[V_T|\mathcal{F}_t] + \mathbb{E}\left[\int_t^T f(U_s,V_s,W_s)\,ds \Big| \mathcal{F}_t\right].$$

**3. Approximation procedure.** In this section we detail the construction of the approximation algorithm of the solution $u$ of $(\mathcal{E})$. We explain how the final form of the discretization procedure can be derived step by step from the forward–backward representation $(E)$. We also present the quantization techniques used in order to compute expectations related to Brownian increments and we discuss the choice of the underlying spatial grids which appear in the approximating scheme.

3.1. *Rough algorithms.*

*Localization procedure.* Recall from the Introduction that the forward–backward equation $(E)$ appears as the starting point of our discretization procedure. Indeed, this couple of stochastic equations provides a probabilistic representation of the quasi-linear PDE $(\mathcal{E})$ and summarizes in an integral form the local evolution of the solution $u$. Define now, for a given integer $N \geq 1$, a regular mesh of $[0,T]$ with step $h \equiv T/N$, that is, set $t_k \equiv kh$, $\forall k \in \{0,\dots,N\}$. Writing the local evolution of $(E)$ and conditioning by $U_{t_k} = x \in \mathbb{R}^d$, we deduce $\forall k \in \{0,\dots,N-1\}$,

(3.1) $U_{t_{k+1}}^{t_k,x} = x + \int_{t_k}^{t_{k+1}} b(U_s^{t_k,x}, V_s^{t_k,x}, W_s^{t_k,x})\,ds + \int_{t_k}^{t_{k+1}} \sigma(U_s^{t_k,x}, V_s^{t_k,x})\,dB_s$

and

$$V_{t_k}^{t_k,x} = \mathbb{E}\left[V_{t_{k+1}}^{t_k,x} + \int_{t_k}^{t_{k+1}} f(U_s^{t_k,x}, V_s^{t_k,x}, W_s^{t_k,x})\,ds\right],$$

$$\mathbb{E}\left[\int_{t_k}^{t_{k+1}} W_s^{t_k,x}\,ds\right] = \mathbb{E}[V_{t_{k+1}}^{t_k,x}(B_{t_{k+1}} - B_{t_k})] + O(h^{3/2}),$$

where the superscript $(t_k,x)$ denotes the starting point of the diffusion process $U$. The remaining term $O(h^{3/2})$ is a consequence of Assumption (A), (2.1) (relationships between $V, W$ and $u$) and Theorem 2.1 (boundedness of $u$ and $\nabla_x u$). Relation (2.1) also yields

$$u(t_k,x) = \mathbb{E}\left[u(t_{k+1}, U_{t_{k+1}}^{t_k,x}) + \int_{t_k}^{t_{k+1}} f(U_s^{t_k,x}, V_s^{t_k,x}, W_s^{t_k,x})\,ds\right],$$

(3.2)
$$\mathbb{E}\left[\int_{t_k}^{t_{k+1}} W_s^{t_k,x}\,ds\right] = \mathbb{E}[u(t_{k+1}, U_{t_{k+1}}^{t_k,x})(B_{t_{k+1}} - B_{t_k})] + O(h^{3/2}).$$

In the following the Brownian increment $B_{t_{k+1}} - B_{t_k}$ is denoted by $\Delta B^k$. In particular, we derive from the above relation that, neglecting the rest, the best constant approximation of $(W_s^{t_k,x})_{s \in [t_k, t_{k+1}]}$ in the $L^2([t_k, t_{k+1}] \times \Omega, ds \otimes d\mathbb{P})$ sense is given by

$$(3.3) \qquad \hat{W}_{t_k}^{t_k,x} \equiv h^{-1} \mathbb{E}[u(t_{k+1}, U_{t_{k+1}}^{t_k,x}) \Delta B^k].$$

Relationships (3.1), (3.2) and (3.3) provide a rough background to discretize the FBSDE ($E$). However, this first form is not satisfactory from an algorithmic point of view. Indeed, because of the strong coupling between the forward and the backward equations, the transition of the diffusion depends on the solution itself, both in the drift term and in the martingale part. At the opposite, in the so-called "pure backward" case, or, correspondingly, for semi-linear equations, the underlying operator does not depend on the solution. In such a case, the classical Euler machinery applies to discretize the decoupled diffusion $U$.

*Induction principle.* Recall that similar difficulties occur to establish the unique solvability of the FBSDE ($E$). In [6] the first author overcomes the strong coupling between the forward and backward equations by solving by induction the local versions of ($E$) on $[t_k, t_{k+1}]$, $k$ running downward from $N-1$ to 0. By analogy with this approach, the discretization procedure of the forward component on a step $[t_k, t_{k+1}]$, $0 \le k \le N-1$, must take into account the issues of the former local discretizations of the backward equation and, more specifically, the approximations of $u(t_{k+1}, \cdot)$ and $v(t_{k+1}, \cdot)$.

*Predictors.* Assume to this end that, at time $t_{k+1}$, some approximations $\bar{u}(t_{k+1}, \cdot), \bar{v}(t_{k+1}, \cdot)$ of $u(t_{k+1}, \cdot)$, $v(t_{k+1}, \cdot)$ are available on the whole space. These approximations appear as the "natural" predictors of the true solution and of its gradient on $[t_k, t_{k+1}]$. Introducing the forward approximating transition

$$(3.4) \quad \mathcal{T}(t_k, x) \equiv b(x, \bar{u}(t_{k+1}, x), \bar{v}(t_{k+1}, x))h + \sigma(x, \bar{u}(t_{k+1}, x))\Delta B^k,$$

we derive an associated updating procedure by setting

$$(3.5) \quad \begin{aligned} \bar{u}(t_k, x) &\equiv \mathbb{E}[\bar{u}(t_{k+1}, x + \mathcal{T}(t_k, x))] + hf(x, \bar{u}(t_{k+1}, x), \bar{v}(t_{k+1}, x)), \\ \bar{v}(t_k, x) &\equiv h^{-1} \mathbb{E}[\bar{u}(t_{k+1}, x + \mathcal{T}(t_k, x))\Delta B^k]. \end{aligned}$$

Once the predictors are updated, the procedure can be iterated. Of course, at time $T = t_N$, we set $\bar{u}(t_N, \cdot) \equiv H(\cdot)$ and $\bar{v}(t_N, \cdot) \equiv \nabla_x H(\cdot)\sigma(\cdot, H(\cdot))$. Note, in particular, that the expectations appearing in (3.5) are correctly defined. Indeed, a simple induction procedure shows from Assumption (A) that $\bar{u}$ and $\bar{v}$ are bounded on $\{t_0, \ldots, t_N\} \times \mathbb{R}^d$ (but the bound depends on the discretization parameters).

*Spatial discretization.* To obtain a numerical scheme, the most natural strategy consists in defining the approximations $\bar{u}(t_k, \cdot)$ and $\bar{v}(t_k, \cdot)$ of the true solution and its gradient on a discrete subset of $\mathbb{R}^d$. Those approximations could then be extended to the whole space with a linear interpolation procedure. However, in high dimension, this last operation can be computationally demanding. We thus prefer, for simplicity, to restrict the approximations to a given spatial grid $\mathcal{C}_k \equiv \{(x_j^k)_{j \in \mathcal{I}_k}, \mathcal{I}_k \subset \mathbb{N}^*\} \subset \mathbb{R}^d$, for $k \in \{0, \dots, N\}$. This choice imposes to modify (3.5). Indeed, the "terminal" value $x + \mathcal{T}(t_k, x)$ must belong to the former grid $\mathcal{C}_{k+1}$.

Hence, denoting by $\Pi_{k+1}$ a projection mapping on the grid $\mathcal{C}_{k+1}$, we replace (3.5) by, $\forall x \in \mathcal{C}_k$,

$$
\begin{aligned}
(3.6) \quad \bar{u}(t_k, x) &\equiv \mathbb{E}[\bar{u}(t_{k+1}, \Pi_{k+1}(x + \mathcal{T}(t_k, x)))] \\
&\quad + h f(x, \bar{u}(t_{k+1}, x), \bar{v}(t_{k+1}, x)), \\
\bar{v}(t_k, x) &\equiv h^{-1} \mathbb{E}[\bar{u}(t_{k+1}, \Pi_{k+1}(x + \mathcal{T}(t_k, x))) \Delta B^k].
\end{aligned}
$$

In the following, we suppose that $\forall (i, j) \in \{0, \dots, N\}^2$, $j < i \Rightarrow \mathcal{C}_j \subset \mathcal{C}_i$, so that $\bar{u}(t_{k+1}, x), \bar{v}(t_{k+1}, x)$ are well defined for $x \in \mathcal{C}_k$. Note that, if the cardinal of $\mathcal{C}_k$ is finite for every $k$, the above scheme is already implementable up to the computations of the underlying expectations.

*Global updating.* The use of the predictors $\bar{u}(t_{k+1}, \cdot), \bar{v}(t_{k+1}, \cdot)$ is an alternative to the standard fixed point procedure. This latter consists in giving first some global predictors $\bar{u}^0(t_k, \cdot), \bar{v}^0(t_k, \cdot)$, $k \in \{0, \dots, N\}$. These are used to compute the transitions of the approximating forward process. In this way, we obtain a decoupled forward–backward system, whose solution may be computed by a standard dynamic programming algorithm. A complete descent of this algorithm from $k = N$ to $k = 0$ produces $\bar{u}^1(t_k, \cdot), \bar{v}^1(t_k, \cdot)$, $k \in \{0, \dots, N\}$, from which we can iterate the previous procedure. In this frame, the underlying distance used to describe the convergence of the fixed point procedure involves all the discretization times and all the spatial points. This strategy appears as a "global updating" one.

From a numerical point of view, this seems unrealistic. Indeed, one would need to solve a large number of linear problems. This would either require to use massive Monte Carlo simulations at each step of the algorithm or to apply, again at each step of the algorithm, a quantization procedure of the approximate diffusion process associated to the current linear problem. Furthermore, it seems intuitively clear that a local updating is far more efficient than a global one.

### 3.2. *Quantization.*

*Expectations approximation.* Two methods are conceivable to compute expectations appearing in (3.6).

The first one consists in applying the classical Monte Carlo procedure for every $k \in \{0, \ldots, N-1\}$ and for every $x \in \mathcal{C}_k$ and, therefore, to repeat this argument $\sum_{k=0}^{N-1} |\mathcal{I}_k|$ times. Such a strategy would lead to perform $\sum_{k=0}^{N-1} |\mathcal{I}_k| \times \varepsilon_{\mathrm{MC}}^{-2}$ elementary operations to compute underlying expectations up to the error term $\varepsilon_{\mathrm{MC}}$. This approach seems rather hopeless.

A more efficient method consists in replacing the Gaussian variables appearing in (3.6) by discrete ones with known weights. This procedure is known as "quantization." Consider to this end a probability measure on $\mathbb{R}^d$ with finite support $(y_i)_{i \in \{1, \ldots, M\}}$ and denote by $(p_i)_{i \in \{1, \ldots, M\}}$ the associated weights. Replace then the Gaussian distribution in (3.6) by this law. For a given $x \in \mathcal{C}_k$, $0 \leq k \leq N$, the expectations appearing in the induction scheme (3.6) then write as computable finite sums.

*Quantization principle.* Generally speaking, for a given random variable $\Delta \in \bigcap_{p \geq 1} L^p(\mathbb{P})$, the quantization procedure consists in replacing $\Delta$ by its projection on a finite grid $\Lambda(M) \equiv \{(y_i)_{i \in \{1, \ldots, M\}}\} \subset \mathbb{R}^d$, $M \in \mathbb{N}^*$. In order to measure the error associated to the grid $\Lambda(M)$, we introduce the so-called "$p$-distortion": $D_{\Delta, p}(\Lambda(M)) \equiv \|\Delta - G_{\Lambda(M)}(\Delta)\|_{L^p(\mathbb{P})}$, $p \geq 1$, where $G_{\Lambda(M)}$ denotes the projection mapping on $\Lambda(M)$. We refer to the monograph of Graf and Luschgy [11] for details.

*Optimal grids.* The crucial step therefore lies in the choice of the grid. The Bucklew–Wise theorem (see Theorem 6.2, Chapter II in [11] for details) then gives, for $\Lambda^*(M)$ achieving the minimum in the $p$-distortion,

$$(3.7) \qquad M^{p/d} D_{\Delta, p}^p(\Lambda^*(M)) \longrightarrow C(p, d) \qquad \text{as } M \to +\infty,$$

where $C(p, d)$ is a constant depending on $p, d$ and the variable at hand.

Various algorithms are available to compute an optimal grid $\Lambda^*(M)$, see, for instance, [2]. We also recall that, for $d > 1$, the optimal grid is not unique.

Up to a rescaling, the basic object associated to Brownian increments is a $d$-dimensional standard normal random variable. Hence, we assume in the following that a grid $\Lambda(M)$ for $\Delta \sim \mathcal{N}(0, \mathbf{I}_d)$, as well as the associated weights $(p_i)_{i \in \{1, \ldots, M\}}$, are given and "perfectly" computed.

*Quantized algorithm.* We are now in position to introduce a more tractable induction principle. Set to this end, for all $k \in \{0, \ldots, N-1\}$, $g(\Delta B^k) \equiv h^{1/2} G_{\Lambda(M)}(h^{-1/2} \Delta B^k)$. Note from the electronic version [9] that, w.l.o.g., for every $p \geq 1$, there exists a constant $C_{\mathrm{Quantiz}}(p, d)$ such that

$$(3.8) \qquad \mathbb{E}[|g(\Delta B^k) - \Delta B^k|^p]^{1/p} \leq C_{\mathrm{Quantiz}}(p, d) h^{1/2} M^{-1/d}.$$

Turn now (3.4) and (3.6) into

$$(3.9) \quad \mathcal{T}(t_k, x) \equiv b(x, \bar{u}(t_{k+1}, x), \bar{v}(t_{k+1}, x)) h + \sigma(x, \bar{u}(t_{k+1}, x)) g(\Delta B^k)$$

and

$$(3.10) \qquad \begin{aligned} \bar{u}(t_k, x) &\equiv \mathbb{E}[\bar{u}(t_{k+1}, \Pi_{k+1}(x + \mathcal{T}(t_k, x)))] \\ &\quad + hf(x, \bar{u}(t_{k+1}, x), \bar{v}(t_{k+1}, x)), \\ \bar{v}(t_k, x) &\equiv h^{-1}\mathbb{E}[\bar{u}(t_{k+1}, \Pi_{k+1}(x + \mathcal{T}(t_k, x)))g(\Delta B^k)]. \end{aligned}$$

To sum up our strategy, the use of predictors allows us to recover a kind of standard dynamic programming principle. The quantization gives an easy, cheap and computable algorithm.

3.3. *Algorithm.* For technical reasons detailed in Section 9, we consider for the convergence analysis a slightly different version of the above algorithm. Namely, we need to change, at a given time $t_k$, the discretization of $b$ and $f$ and, in particular, to replace $\bar{v}(t_{k+1}, \cdot)$ by a new predictor. Concerning the driver of the BSDE, we replace $f(x, \bar{u}(t_{k+1}, x), \bar{v}(t_{k+1}, x))$ by $f(x, \bar{u}(t_{k+1}, x), \bar{v}(t_k, x))$: the definition of $\bar{v}(t_k, x)$ does not involve $\bar{u}(t_k, x)$.

The story is rather different for $b$. Indeed, the definition of $\bar{v}(t_k, x)$ relies on the choice of the underlying transition. In particular, putting $\bar{v}(t_k, x)$ in $b$ as done in $f$ would lead to an implicit scheme.

Nevertheless, for a given intermediate predictor $\hat{v}(t_k, \cdot)$ of $v(t_k, \cdot)$, we can put

$$\mathcal{T}(t_k, x) \equiv b(x, \bar{u}(t_{k+1}, x), \hat{v}(t_k, x))h + \sigma(x, \bar{u}(t_{k+1}, x))g(\Delta B^k).$$

The whole difficulty is then hidden in the choice of $\hat{v}(t_k, x)$. Our strategy consists in choosing $\hat{v}(t_k, x)$ as the expectation of $\bar{v}(t_{k+1}, \cdot)$, with respect to the transition $\mathcal{T}^0(t_k, x) \equiv \sigma(x, \bar{u}(t_{k+1}, x))g(\Delta B^k)$. This transition differs from $\mathcal{T}(t_k, x)$ in the drift $b$ and leads to an explicit scheme. Namely, we set

$$(3.11) \qquad \hat{v}(t_k, x) \equiv \mathbb{E}[\bar{v}(t_{k+1}, \Pi_{k+1}(x + \mathcal{T}^0(t_k, x)))].$$

The predictor $\hat{v}(t_k, \cdot)$ appears as a "regularized" version of $\bar{v}(t_{k+1}, \cdot)$. Thanks to a Gaussian change of variable, the laws of the underlying transitions $\mathcal{T}^0(t_k, x)$ and $\mathcal{T}(t_k, x)$ can be compared; see [9], Section 7.3, for details.

*Final algorithm.*

ALGORITHM 3.1. The final algorithm writes

$$\forall x \in \mathcal{C}_N, \bar{u}(T, x) \equiv H(x), \qquad \bar{v}(T, x) \equiv \nabla_x H(x)\sigma(x, H(x)),$$

$$\forall k \in \{0, \ldots, N-1\}, \ \forall x \in \mathcal{C}_k,$$

$$\mathcal{T}^0(t_k, x) \equiv \sigma(x, \bar{u}(t_{k+1}, x))g(\Delta B^k),$$

$$\hat{v}(t_k, x) \equiv \mathbb{E}[\bar{v}(t_{k+1}, \Pi_{k+1}(x + \mathcal{T}^0(t_k, x)))],$$

$$\mathcal{T}(t_k, x) \equiv b(x, \bar{u}(t_{k+1}, x), \hat{v}(t_k, x))h + \sigma(x, \bar{u}(t_{k+1}, x))g(\Delta B^k),$$

$$\bar{v}(t_k, x) \equiv h^{-1}\mathbb{E}[\bar{u}(t_{k+1}, \Pi_{k+1}(x + \mathcal{T}(t_k, x)))g(\Delta B^k)],$$

$$\bar{u}(t_k, x) \equiv \mathbb{E}[\bar{u}(t_{k+1}, \Pi_{k+1}(x + \mathcal{T}(t_k, x)))] + f(x, \bar{u}(t_{k+1}, x), \bar{v}(t_k, x))h.$$

*A discrete probabilistic representation.* Following the link between ($E$) and ($\mathcal{E}$), define, for $x_0 \in \mathcal{C}_0$, a Markov process on the grids $(\mathcal{C}_k)_{0 \leq k \leq N}$ according to the transitions $(\mathcal{T}(t_k, x))_{k \in \{0, \dots, N-1\}, x \in \mathcal{C}_k}$,

$$(3.12) \quad X_0 \equiv x_0, \ \forall k \in \{0, \dots, N-1\}, \qquad X_{t_{k+1}} \equiv \Pi_{k+1}(X_{t_k} + \mathcal{T}(t_k, X_{t_k})).$$

Referring to the connection between $U$ and $(V, W)$ [see, e.g. (2.1)], put now

$$(3.13) \quad \forall k \in \{0, \dots, N\}, \qquad Y_{t_k} \equiv \bar{u}(t_k, X_{t_k}), \qquad Z_{t_k} \equiv \bar{v}(t_k, X_{t_k}).$$

Note that $Y$ and $Z$ are correctly defined since $X_{t_k}$ belongs to the grid $\mathcal{C}_k$. The couple $(Y, Z)$ appears as a discrete version of the couple $(V, W)$ in ($E$). More precisely, one can prove the following discrete Feynman–Kac formula: $\forall 0 \leq k \leq N - 1$,

$$(3.14) \quad Y_{t_k} = \mathbb{E}\left[H(X_{t_N}) + h \sum_{i=k+1}^{N} f(X_{t_{i-1}}, \bar{u}(t_i, X_{t_{i-1}}), Z_{t_{i-1}}) \Big| \mathcal{F}_{t_k}\right].$$

Note anyhow that the process $Z$ does not appear as the martingale part of the process $Y$. However, thanks to the martingale representation theorem, there exists a progressively measurable process $\overline{Z}$, with finite moment of order two, such that

$$(3.15) \quad Y_{t_N} + h \sum_{i=1}^{N} f(X_{t_{i-1}}, \bar{u}(t_i, X_{t_{i-1}}), Z_{t_{i-1}}) = Y_0 + \int_0^{t_N} \overline{Z}_s \, dB_s.$$

Of course, the process $\overline{Z}$ does not match exactly the process $Z$. However, for a given $k \in \{0, \dots, N-1\}$, it is readily seen from the above expression that the best $\mathcal{F}_{t_k}$-measurable approximation of $(\overline{Z}_s)_{s \in [t_k, t_{k+1}]}$ in $L^2([t_k, t_{k+1}] \times \Omega, ds \otimes d\mathbb{P})$ is given by $h^{-1}\mathbb{E}[Y_{t_{k+1}} \Delta B^k | \mathcal{F}_{t_k}]$. Up to the quantization procedure, this term coincides with $\bar{v}(t_k, X_{t_k})$. In other words, the processes $Z$ and $\overline{Z}$ may be considered as close.

3.4. *Choice of the grids.* Because of the strong coupling, little is a priori known on the behavior of the paths of the forward process. Hence, we cannot compute a kind of optimal grid for $X$. The most natural choice turns out to be the one of Cartesian grids.

*Unbounded Cartesian grids.* Two different choices of grids are conceivable. First, we can treat the case of infinite Cartesian grids: $\forall k \in \{0, \dots, N\}$, $\mathcal{C}_k \equiv \mathcal{C}_\infty$, $\mathcal{C}_\infty \equiv \delta \mathbb{Z}^d$, where $\delta > 0$ denotes a spatial discretization parameter. In this case, the projection mapping writes $\forall x \in \mathbb{R}^d$, $\Pi_\infty(x) \equiv \sum_{y \in \mathcal{C}_\infty} [y \prod_{j=1}^d \mathbf{1}_{[-\delta/2, \delta/2[}(x_j - y_j)]$. In other words, for every $j \in \{1, \dots, d\}$, the coordinate $j$ of $\Pi_\infty(x)$ is given by $(\Pi_\infty(x))_j = \delta \lfloor \delta^{-1} x_j + 1/2 \rfloor$.

This choice actually simplifies the convergence analysis and allows a direct comparison with the results from the existing literature; see [10]. Note, however, that it does not provide a fully implementable scheme since the set $\mathcal{C}_\infty$ is infinite.

*Truncated grids.* Several truncation procedures may be considered, but all need to take into account the specific geometry of a nondegenerate diffusion, or, more simply, of the Brownian motion. Set, for example, for a given $R > 0$, and for all $i \in \{0, \dots, N\}$, $\mathcal{C}_i \equiv \mathcal{C}_\infty \cap \Delta_i$, where

$$
\begin{aligned}
(3.16) \quad & \Delta_i \equiv \{x \in \mathbb{R}^d, \forall 1 \le j \le d, \\
& -\delta \lfloor (R + \rho \psi(t_i)) \delta^{-1} \rfloor - \delta/2 \le x_j < \delta \lfloor (R + \rho \psi(t_i)) \delta^{-1} \rfloor + \delta/2\},
\end{aligned}
$$

where $\psi(t) = t^\eta \mathbf{1}_{\{t > 0\}}$, $\eta \in [0, 1/2)$, is meant to take into account the Hölder regularity of the Brownian path. The larger is $\eta$, the smaller is the number of points involved in the discretization procedure. However, since the proof of the convergence of the algorithm is far from being trivial, we restrict our analysis to the case $\eta = 0$.

Note also that the particular choice of the bounds in the definition of $\Delta_i$ ensures that for all $x \in \mathbb{R}^d$, $\Pi_\infty(x) \in \mathcal{C}_i \Leftrightarrow x \in \Delta_i$. Hence, for every $i \in \{0, \dots, N\}$, $\Pi_i$ writes

$$
\begin{aligned}
(3.17) \quad & \forall 0 \le i \le N, \ \forall x \in \Delta_i, & & \Pi_i(x) \equiv \mathcal{Q}(R + \rho, \Pi_\infty(x)) = \Pi_\infty(x), \\
& \forall 1 \le i \le N, \ \forall x \notin \Delta_i, & & \Pi_i(x) \equiv \mathcal{Q}(R + \rho, \Pi_\infty(x)), \\
& \forall x \notin \Delta_0, & & \Pi_0(x) \equiv \mathcal{Q}(R, \Pi_\infty(x)),
\end{aligned}
$$

where, for a given $(r, y) \in \mathbb{R}^{+*} \times \mathbb{R}^d$, $\mathcal{Q}(r, y)$ denotes the orthogonal projection of $y$ on the hypercube $[-\delta \lfloor r \delta^{-1} \rfloor, \delta \lfloor r \delta^{-1} \rfloor]^d$: $\mathcal{Q}(r, y) \equiv ((y_i \vee (-\delta \lfloor r \delta^{-1} \rfloor)) \wedge (\delta \lfloor r \delta^{-1} \rfloor))_{1 \le i \le d}$. Note finally that $R$ is fixed by the reader once for all in function of the set on which $u$ has to be approximated at the initial time. At the opposite, $\rho$ appears as a discretization parameter chosen by the reader in function of the required precision and of the affordable complexity for Algorithm 3.1.

## 4. Convergence results.
This section is devoted to the convergence analysis of $\bar{u}$ to $u$. As stated in the following theorem, which is the main result of the paper, five different types of errors can be distinguished:

THEOREM 4.1. *Let $p \geq 2$. There exist two constants $c_{4.1}$ and $C_{4.1}$, only depending on $p$, $T$ and on known parameters deriving from Assumption (A), such that, for $h < c_{4.1}$, $\delta^2 < h$, $M^{-2/d} < h$ and $\rho \geq 1$,*

$$\sup_{x \in \mathcal{C}_0} |u(0, x) - \bar{u}(0, x)|^2 \leq C_{4.1} \mathcal{E}^2 (\text{global}),$$

*with $\mathcal{E}^2(\text{global}) \equiv \mathcal{E}^2(\text{time}) + \mathcal{E}^2(\text{space}) + \mathcal{E}^2(\text{trunc}) + \mathcal{E}^2(\text{quantiz}) + \mathcal{E}^2(\text{gradient}, p)$ and $\mathcal{E}(\text{time}) \equiv h^{1/2}$, $\mathcal{E}(\text{space}) \equiv h^{-1}\delta$, $\mathcal{E}(\text{trunc}) \equiv R/(R+\rho)$, $\mathcal{E}(\text{quantiz}) \equiv h^{-1/2}M^{-1/d}$, $\mathcal{E}(\text{gradient}, p) \equiv h^{p/2 + d/4 - 1/2} M^{-p/d} \delta^{-p - d/2}$.*

REMARK 4.1. The FBSDE counterpart of Theorem 4.1 is given in Section 4.3: see Theorems 4.2 and 4.3.

4.1. *Classification of errors.* We now detail the meaning of the different errors appearing in Theorem 4.1:

*Temporal discretization error $\mathcal{E}(\text{time})$.* The $1/2$ exponent appearing in the definition of $\mathcal{E}(\text{time})$ corresponds to the Hölder regularity of $u$ and $\nabla_x u$ in time and to the $L^2(\mathbb{P})$ $1/2$-Hölder property of the Brownian increments.

*Spatial discretization error $\mathcal{E}(\text{space})$.* This quantity highly depends on the ratio between the spatial and the temporal steps. This connection between $\delta$ and $h$ can be explained as follows: the drift part of the transitions $(\mathcal{T}(t_k, \cdot))_{0 \leq k \leq N}$ is of order $h$ and the diffusive one is of order $h^{1/2}$. Thus, to take into account the influence of the drift at the local level, the spatial discretization parameter must be smaller than $h$. In other words, $\delta h^{-1}$ must be small.

*Quantization error $\mathcal{E}(\text{quantiz})$.* This error depends on the ratio between the distortion and the temporal step. The quantity $\mathcal{E}(\text{quantiz})$ represents the typical bound between $\bar{v}(t_k, X_{t_k})$ and the best $\mathcal{F}_{t_k}$ measurable approximation of the process $(\overline{Z}_s)_{s \in [t_k, t_{k+1}]}$, that is, between $\bar{v}(t_k, X_{t_k})$ and $h^{-1}\mathbb{E}[\bar{u}(t_{k+1}, \Pi_{k+1}(X_{t_k} + \mathcal{T}(t_k, X_{t_k})))\Delta B^k | \mathcal{F}_{t_k}]$. Note, indeed, that the distance between $\Delta B^k$ and $g(\Delta B^k)$ is of order $h^{1/2}M^{-1/d}$, see (3.8). Since the underlying expectation is divided by $h$, this leads to a term in $h^{-1/2}M^{-1/d}$.

*Truncation error $\mathcal{E}(\text{trunc})$.* As written in Theorem 4.1, it depends on $R$ and $\rho$, where $R$ denotes the radius of the initial grid $\mathcal{C}_0$ and $R + \rho$ the radius of the grids $(\mathcal{C}_k)_{1 \leq k \leq N}$. If $\rho$ tends to $+\infty$, that is, if the grids are not truncated, this error term reduces to zero.

Generally speaking, $\mathcal{E}(\text{trunc})$ appears as the Bienaymé–Chebyshev estimate of the probability that the approximating process $X$ stays inside the grids $(\mathcal{C}_k)_{0 \leq k \leq N}$. The lack of relevant estimates of the discretized version of the drift $b$ (recall that the function $b$ is not bounded) and, more specially, of the discretized gradient $\bar{v}$, explains the reason why the Bienaymé–Chebyshev estimate applies in this framework and not better ones (as the Bernstein inequality). We also recall that the unboundedness of the coefficients is the most common case in the applications, see, for example, Section 5.2.

*Gradient error $\mathcal{E}$(gradient, $p$).*  This extra error is generated by the lack of estimates of the discretized gradient $\bar{v}$. This term follows from the specific choice of the predictor $\hat{v}$ made in Section 3.3 and appears in the second step of the proof of Theorem 4.1; see, more precisely, Sections 7.1 and 7.3.

The convergence of $\mathcal{E}$(gradient, $p$) toward 0 relies on the term $h^{p/2} M^{-p/d} \delta^{-p}$, $M$ being chosen large enough and $p$ as large as necessary. In short, this reduced form represents the probability that the distance between the Gaussian increment and its quantization exceeds the spatial step $\delta$. Note, indeed, from (3.8) that, for every $p \geq 2$, $\mathbb{P}\{|\Delta B^k - g(\Delta B^k)| > \delta\} \leq C_{\text{Quantiz}}(p, d) h^{p/2} M^{-p/d} \delta^{-p}$. Thus, the error term $\mathcal{E}$(gradient, $p$) depends on the ratio between the spatial discretization step and the quantization distortion of the underlying Gaussian increments.

The above probability appears in the control of the distance between the predictor $\hat{v}$ and the true gradient $v$. In this frame, the strategy consists in writing the predictor $\hat{v}$ as an expectation with respect to the Gaussian kernel and not to its quantized version. Generally speaking, this strategy holds when the quantized transition $\mathcal{T}(t_k, x)$ and its Gaussian counterpart belong to the same cell of the spatial grid, that is, when the distance between the Brownian increment and the quantized one is of the same order as the length of a given cell. Since the spatial grid step is given by $\delta$, we then need to control the probability that the difference between the increments exceeds $\delta$.

Of course, when $b$ does not depend on $z$, there is no reason to define $\hat{v}$. In such a case, $\mathcal{E}$(gradient, $p$) reduces to 0.

### 4.2. *Comments on the rate of convergence.*

*Error in function of $h$.*  To detail in a more explicit way the rate of convergence given by Theorem 4.1, we give an example in which $\rho$ ($\rho < +\infty$), $\delta$ and $M$ are expressed as powers of $h$. Assume, indeed, that $\rho$, $\delta$ and $M$ are chosen in the following way: $\rho = R h^{-1/2}$, $\delta \equiv h^{1+\gamma}$, $M^{-2/d} \equiv h^{1+\beta}$, $\gamma, \beta \geq 0$.

In such a case, $\mathcal{E}$(gradient, $p$) = $\exp[\ln(h)[p(\beta/2 - \gamma) - (d/2 + 1 + \gamma d)/2]]$. To ensure the convergence of the algorithm, we then need to choose

$$p(\beta/2 - \gamma) - (d/2 + 1 + \gamma d)/2 > 0 \quad \Longleftrightarrow \quad \beta > 2\gamma + (1/p)(d/2 + 1 + \gamma d).$$

Put finally $\beta = 2\gamma + (d/2 + 1 + \gamma d)/p + \eta$, $\eta > 0$. The rate of convergence of the fully implementable algorithm is given by $\sup_{x \in \mathcal{C}_0} |u(0,x) - \bar{u}(0,x)|^2 \leq C_{4.1}[h + h^{2\gamma} + h^\beta + h^{p\eta}]$.

Taking $\gamma = 1/2$ and $\eta = 1/p$ then yields $\sup_{x \in \mathcal{C}_0} |u(0,x) - \bar{u}(0,x)|^2 \leq C_{4.1} h$. In particular, for $p$ large enough, the exponent $\beta$ is close to 1 and the number $M$ of points needed to quantify the Brownian increments is close to $h^{-d}$. Here is the limit of the method: for a large $d$ and a small $h$, we need a rather large number of points for the Gaussian quantization. Recall anyhow that the Gaussian grids are computed once for all. Thus, the numerical effort to get sharp quantization grids can be made apart from our algorithm.

*Estimates of $\nabla_x u$.* The reader might wonder about the estimate of the gradient of $u$. Note in this framework that two strategies are conceivable.

First, the probabilistic counterpart of Theorem 4.1 given in Section 4.3 provides an $L^2$ estimate of the distance between $\bar{v}$ and the gradient of the true solution. Note, however, that the underlying $L^2$ norm is taken with respect to the distribution of the discrete process $X$ [cf. (3.12)].

To get a joint estimate of the solution and of its gradient with respect to the *supremum* norm, the reader can apply the following strategy: differentiate if possible the PDE ($\mathcal{E}$) and apply, once again if possible, Algorithm 3.1 to $(u, \nabla_x u)$, seen as the solution of a system of parabolic quasi-linear PDEs. Such a strategy is applied in Section 5 to the solution of the porous media equation and to its gradient. Note that this approach coincides with the one followed by Douglas, Ma and Protter [10].

4.3. *Estimates of the discrete processes.* We now translate Theorem 4.1 in a more probabilistic way. Recall indeed that, in several situations (e.g., in financial mathematics), the knowledge of the triple $(U, V, W)$ is as crucial as the knowledge of the couple $(u, \nabla_x u)$.

We then prove that $(X, Y, Z)$ and $(U, V, W)$ get closer in a suitable sense as $h$, $\delta$, $M^{-1}$ and $\rho^{-1}$ vanish. Note, however, that we are not able to prove that the distance between $(X, Y, Z)$ and $(U, V, W)$ over the whole interval $[0,T]$ tends to zero. Indeed, since the projections $(\Pi_i)_{0 \leq i \leq N}$ map every point outside the sets $(\Delta_i)_{0 \leq i \leq N}$ onto the boundaries of $(\mathcal{C}_i)_{0 \leq i \leq N}$ [see, e.g., (3.17)], we do not control efficiently the transition of the process $X$ after the first hitting time of the boundaries of the grids by $X$. It is then well understood that we have to stop the triple $(X, Y, Z)$ at this first hitting time. Put to this end

$$\tau_\infty \equiv \inf\{(t_k)_{1 \leq k \leq N}, X_{t_{k-1}} + \mathcal{T}(t_{k-1}, X_{t_{k-1}}) \notin \Delta_k\}, \qquad \inf(\varnothing) = +\infty.$$
(4.1)

First, as a bypass product of the proof of Theorem 4.1, the function $\bar{v}$ provides an approximation of $v$ in the following $L^2$ sense:

THEOREM 4.2. *Let $p \geq 2$. Then, there exist two constants $c_{4.2}$ and $C_{4.2}$, only depending on $p$ and on known parameters deriving from Assumption* (A), *such that, for $h < c_{4.2}$, $\delta^2 < h$, $M^{-2/d} < h$ and $\rho \geq 1$,*

$$h \sum_{i=0}^{N-1} \mathbb{E}[|\bar{v}(t_i, X_{t_i})\mathbf{1}_{\{t_i < \tau_\infty\}} - v(t_i, X_{t_i})|^2] \leq C_{4.2}\mathcal{E}^2(\text{global}).$$

Moreover, the triple $(X, Y, Z)$ stopped at time $\tau_\infty$ satisfies the following:

THEOREM 4.3. *Let $p \geq 2$. Then, there exist two constants $c_{4.3}$ and $C_{4.3}$, only depending on $p$ and on known parameters deriving from Assumption* (A), *such that, for $h, \delta, M$ as in the previous theorem*

$$\mathbb{E}\left[\sup_{i \in \{0,...,N\}} |X_{t_i \wedge \tau_\infty} - U_{t_i}|^2\right] + \mathbb{E}\left[\sup_{i \in \{0,...,N\}} |Y_{t_i \wedge \tau_\infty} - V_{t_i}|^2\right]$$

$$+ h \sum_{i=0}^{N-1} \mathbb{E}[|Z_{t_i}\mathbf{1}_{\{t_i < \tau_\infty\}} - W_{t_i}|^2] \leq C_{4.3}\mathcal{E}^2(\text{global}).$$

**5. Numerical examples.** In order to compare the results we obtain with our algorithm to a reference value, we choose equations that admit an explicit solution. In this frame, we focus on three examples: the one-dimensional Burgers equation, the deterministic KPZ equation in dimension two and the one-dimensional porous media equation.

5.1. *One-dimensional Burgers equation.* Consider first the backward Burgers equation:

$$\partial_t u(t,x) - (u\partial_x u)(t,x) + \frac{\varepsilon^2}{2}\partial^2_{x,x}u(t,x) = 0,$$

(5.1)
$$(t,x) \in [0,T[\times \mathbb{R}, \varepsilon > 0$$

$$u(T,x) = H(x), \qquad x \in \mathbb{R}, H \in C_b^{2+\alpha}(\mathbb{R}), \alpha \in ]0,1[.$$

Using a nonlinear transformation, one can derive an explicit expression of the solution of (5.1). This is known as the Cole–Hopf factorization, see [23], Chapter IV, or [24], Chapter III, for details. The solution of (5.1) then writes

(5.2) $\quad \forall (t,x) \in [0,T] \times \mathbb{R}, \qquad u(t,x) = \dfrac{\mathbb{E}[H(x + \varepsilon B_{T-t})\phi(x + \varepsilon B_{T-t})]}{\mathbb{E}[\phi(x + \varepsilon B_{T-t})]},$

where $B$ is a standard Brownian motion and

$$\forall y \in \mathbb{R}, \qquad \phi(y) \equiv \exp\left(-\varepsilon^{-2}\int_0^y H(u)\,du\right).$$

From the explicit representation (5.2), we can derive numerically, using, for example, a Riemann sum, a Monte Carlo method or a quantized version of the expectation (5.2), a reference solution to test the algorithm.

The reader may object that the Burgers equation is actually semi-linear and not quasi-linear. Actually, it depends on whether we consider the nonlinear term as a drift or as a second member. We describe below the algorithms associated to these two points of view, even if the coupled case is the only one to fulfill Assumption (A).

Moreover, in the forward–backward representation of the Burgers equation, the estimation procedure of the gradient is not necessary to compute the approximate solution $\bar{u}$. Numerically, this case turns out to be the most robust. Finally, in both cases, the intermediate predictor $\hat{v}$ is useless: in the coupled case, the drift of the diffusion $U$ reduces to $V$ (and thus does not depend on $W$), and in the decoupled one, the drift vanishes.

5.1.1. *Explicit expression of the algorithms.* For a given final condition $H \in C_b^{2+\alpha}(\mathbb{R})$, $\alpha \in ]0,1[$, we write the following:

ALGORITHM 5.1 (Coupled case).

$\forall x \in \mathcal{C}_N, \ \bar{u}(T,x) \equiv H(x)$,

$\forall k \in \{0,\dots,N-1\}, \ \forall x \in \mathcal{C}_k$,

$\quad \bar{u}(t_k,x) \equiv \mathbb{E}[\bar{u}(t_{k+1},\Pi_{k+1}(x - \bar{u}(t_{k+1},x)h + \varepsilon g(\Delta B^k)))]$,

$\quad \bar{v}(t_k,x) \equiv h^{-1}\mathbb{E}[\bar{u}(t_{k+1},\Pi_{k+1}(x - \bar{u}(t_{k+1},x)h + \varepsilon g(\Delta B^k)))g(\Delta B^k)]$.

ALGORITHM 5.2 (Pure backward case).

$\forall x \in \mathcal{C}_N, \ \bar{u}(T,x) \equiv H(x)$,

$\forall k \in \{0,\dots,N-1\}, \ \forall x \in \mathcal{C}_k$,

$\quad \bar{u}(t_k,x) \equiv \mathbb{E}[\bar{u}(t_{k+1},\Pi_{k+1}(x + \varepsilon g(\Delta B^k)))] - h\varepsilon^{-1}\bar{u}(t_{k+1},x)\bar{v}(t_k,x)$,

$\quad \bar{v}(t_k,x) \equiv h^{-1}\mathbb{E}[\bar{u}(t_{k+1},\Pi_{k+1}(x + \varepsilon g(\Delta B^k)))g(\Delta B^k)]$.

5.1.2. *Numerical results.* In order to avoid first to truncate the grids, we choose a periodic initial solution. Put to this end $H(x) = \sin(2\pi x)$ and derive from (5.2) that $u$ is 1-periodic. This allows to define $\bar{u}(t_k,\cdot)$ on $\mathcal{C}_\infty$ by setting $\forall x \in \mathcal{C}_\infty, \ \bar{u}(t_k,x) \equiv \bar{u}(t_k,x - \lfloor x \rfloor)$. Hence, we can set $\mathcal{C}_k \equiv \mathcal{C}_\infty$ for $k \in \{0,\dots,N-1\}$. For $T = 1, \delta = 10^{-3}, \ h = 0.01, \ M = 160, \ \varepsilon = 0.15$, we present below the results of the previous algorithms. The explicit solution given by (5.2) is approximated by quantization techniques with a 500 points grid. We plot below some profiles of the reference value for various discretization times, as well as the pointwise absolute error between this reference solution and the approximations obtained with our algorithms. See Figure 1.

On the profiles of the explicit solution, the abscises of the peaks of the initial sinusoidal wave are going closer to each other up to a given time $t_0$. This is a typical shocking wave behavior. Because of the viscosity, that is, $\varepsilon$ is nonzero, there is no shock and the amplitude of the wave decays when $t$ goes to zero.

From a numerical point of view, the coupled case provides several advantages. First, the convergence of Algorithm 5.1 does not rely on the discretization procedure of the gradient. In short, there is no reason to update the gradient in order to obtain the approximate solution with the first algorithm. The computation of $\bar{v}$ just provides in this case an $L^2$ estimate of the gradient. At the opposite, this computation is necessary in Algorithm 5.2.

Moreover, since the coefficient $f(y, z) = \varepsilon^{-1} yz$ is not globally Lipschitz in the pure backward case, it is then another story to establish the convergence of Algorithm 5.2.

These theoretical remarks are confirmed by the pictures below. Even though Algorithm 5.2 does not behave too poorly, it is still less precise than Algorithm 5.1. The factor between the absolute pointwise errors of the two algorithms is approximately 5.

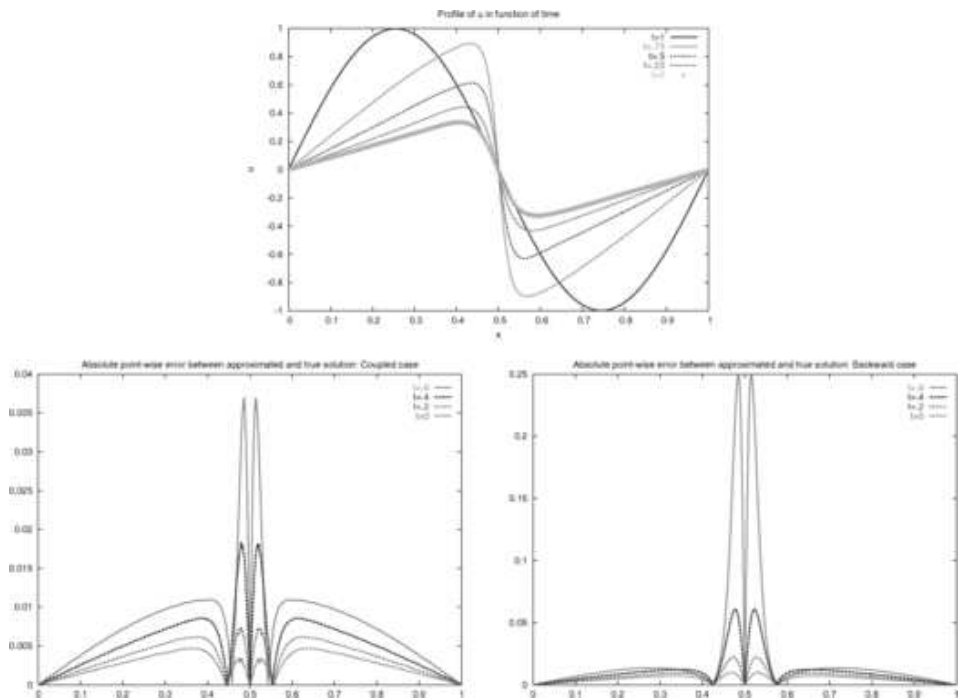

Fig. 1.

*Truncation error.* We now illustrate the effects of truncation and deal with a nonperiodic final data. Namely, we take $H(x) = \exp(-x^2/2)$, $T = 1$, $h = 0.02$, $\rho = 3$, $\delta = \rho/500$, $M = 250$. The reference value, see profiles below, is computed from the Cole–Hopf explicit solution by quantization techniques with a 500 points grid. We run Algorithm 3.1 with the previous parameters to obtain Figure 2.

Choose now $R = 1$: the expected truncated error $\mathcal{E}(\text{trunc})$ is given by 0.25, whereas the absolute point-wise error between both solutions is bounded by 0.05 on $[-1, 1]$. This emphasizes the difficulty to control the truncation procedure in our algorithm. There are two possible arguments to explain this difference between 0.25 and 0.05. First, as explained in Section 9.1, our way to estimate $\mathcal{E}(\text{trunc})$ is suitable for unbounded drifts $b$ and, more particulary, for drifts depending on the gradient. In our case, the drift is bounded (since the solution is bounded by 1), and most relevant estimates could apply. Second, the fast decay of the final condition $H$ may explain the low influence of distant points on the values of the solution on $[-1, 1]$.

Note also that the relative error is close to 0.1 on $[-1, 1]$. A possible strategy to decrease it would consist in refining the spatial mesh.

We also feel that the choice of the rough projection mappings $(\Pi_k)_{0 \le k \le N}$ deeply affects the global error. To investigate more precisely their influence, we replace them by standard linear interpolation procedures (which are defined in an obvious way since the underlying space is one dimensional). In short, this permits to extend continuously the approximated solution $\bar{u}$ to the whole space. With the same parameters as above, we then get Figure 3.

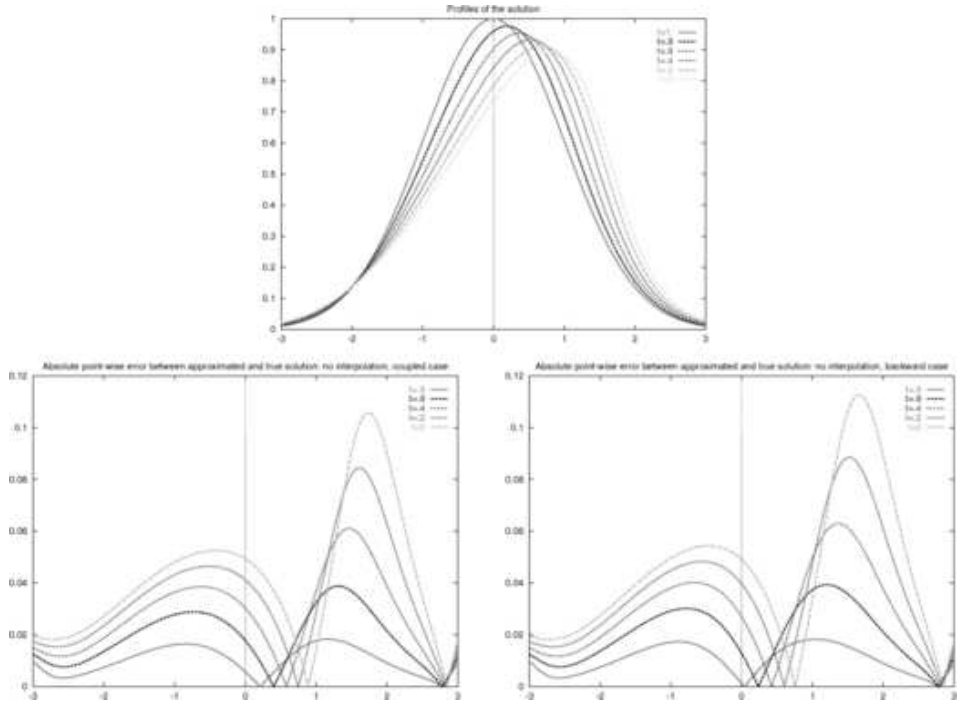

Fig. 2.

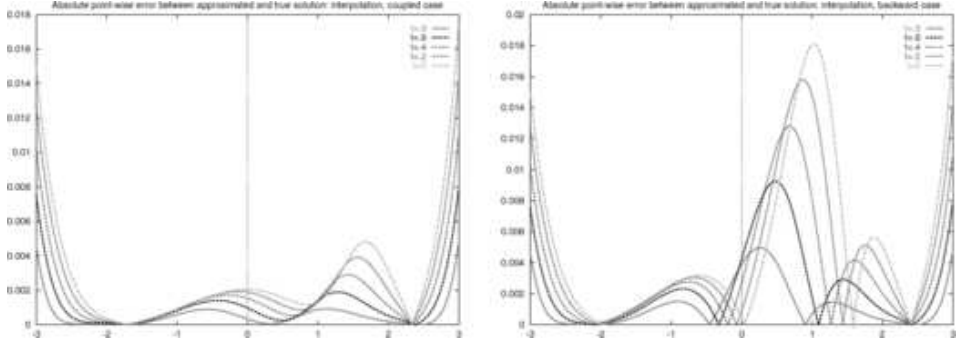

Fig. 3.

Numerically, the interpolation can thus be really relevant to improve the convergence (see Section 9.2 for further details and explanations on this point). To obtain the same precision without interpolation, we need to refine significantly the parameters (taking, e.g., $\delta = 2 \times 10^{-4}$). Let us finally mention that the results obtained with the coupled representation and the linear interpolation are still more accurate than with the backward one.

5.2. *Deterministic KPZ equation.* In this subsection we focus on the so-called "deterministic KPZ" equation (see, e.g., [12] and [24], Chapter I, for a physical interpretation):

$$\partial_t u(t,x) + \frac{1}{2}\mathrm{tr}(\sigma\sigma^*\nabla^2_{x,x}u(t,x)) + \frac{\nu}{2}|\sigma^*\nabla_x u(t,x)|^2 = 0,$$

(5.3)
$$(t,x) \in [0,T[\times\mathbb{R}^d,$$

$$u(T,x) = H(x), \qquad x \in \mathbb{R}^d,$$

where $\nu \in \mathbb{R}^{+*}$ is a given parameter and $\sigma$ a given constant matrix such that $\sigma\sigma^*$ is positive definite.

Such an equation admits too a "Cole–Hopf explicit solution" (see again [12]) that writes $u(t,x) = \nu^{-1}\log(\mathbb{E}[\exp(\nu H(x + \sigma B_{T-t}))])$. We then apply Algorithm 3.1 to (5.3) seen as a true quasi-linear equation (so-called "coupled case" in the former subsection).

Concerning the initial condition, we choose $H(x) = \prod_{i=1}^{d}\sin(2\pi x_i)$. By construction, we have $\forall x \in \mathbb{R}^d, \ \forall k \in \mathbb{Z}^d, \ u(t, x+k) = u(t,x)$. Since the solution is periodic, $\bar{u}$ can be defined on the whole grid $\mathcal{C}_\infty$ (see also Section 5.1.2). We now present the results for $d = 2$, $\nu = 0.3$, $T = 0.5$, $h = 0.02$, $\delta = 5 \times 10^{-4}$, $M = 160$ and $\sigma\sigma^* = \begin{pmatrix} 1 & \theta \\ \theta & 1 \end{pmatrix}$ with $\theta = 0.8$. The reference value and its gradient have been derived from the explicit writing of $u$ using quantization techniques with a 500 points grid. At $t = 0$, one has Figure 4.

The relative error between the approximate and true solutions is at most 0.25. The explanation seems rather simple: the explicit solution quickly decays as time decreases. Anyway, we feel that our algorithm manages to catch this specific decreasing phenomenon.

Let us also mention that the last picture represents the pointwise difference of the true and approximated gradients, but the control given by Theorem 4.2 just holds in $L^2$.

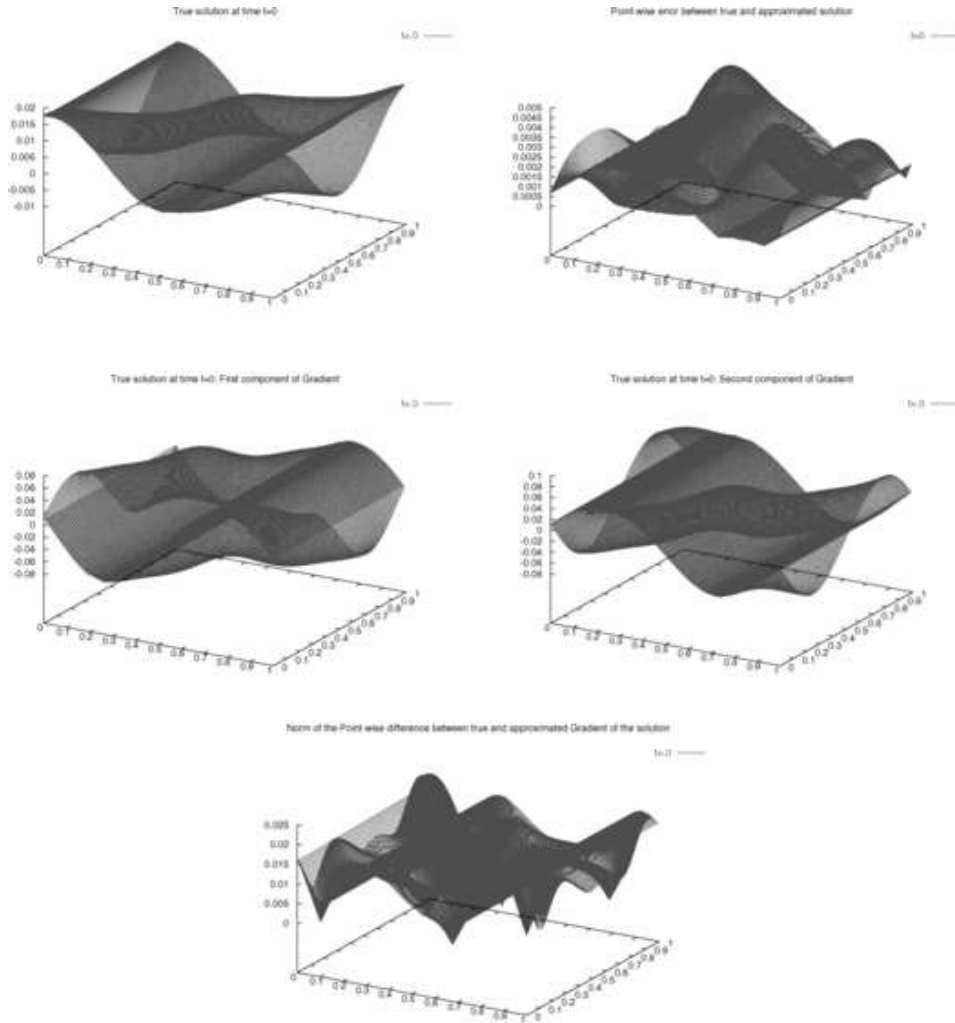

Fig. 4.

5.3. *Porous media equation.* To conclude this section, we focus on the equation (this example is taken from [16])

$$\partial_t u(t,x) + (u \partial_{x,x}^2 u)(t,x) + (\partial_x u)^2(t,x) + u^2(t,x) = 0,$$

$$(5.4) \qquad\qquad (t,x) \in ]0,T] \times \mathbb{R},$$

$$u(T,x) = T^{-1} \frac{4}{3} \cos^2\left(\frac{\pi x}{L}\right), \qquad L = 2\sqrt{2}\pi,$$

which admits the $L$-periodic explicit solution

$$u(t,x) = t^{-1} \frac{4}{3} \cos^2\left(\frac{\pi x}{L}\right).$$

Note that (5.4) does not fulfill Assumption (A). In the sequel, we choose without any rigorous justifications to apply Algorithm 3.1 on $[T/2,T]$ (note, however, for a rough explanation that the quadratic growth of the coefficients ensures that Theorem 4.1 holds on a suitable interval $[t,T]$, for $t$ close enough to $T$ and, in the same way, Theorem 2.1 applies away from 0).

Nevertheless, as explained in Section 4.2, this procedure just provides an $L^2$-estimate of $\nabla_x u$. In this framework, we have decided to apply the so-called "differentiated" approach, described in Section 4.2, to obtain a pointwise estimate of $\nabla_x u$ (see Algorithm 5.3 below).

Note finally from the periodicity of $u$ that $\bar{u}$ can be defined on the whole grid $\mathcal{C}_\infty$ as in the previous example (see also Section 5.1.2).

ALGORITHM 5.3 (Differentiated algorithm).

$$\forall x \in \mathcal{C}_N, \qquad \bar{u}(T,x) = T^{-1} \frac{4}{3} \cos^2\left(\frac{\pi x}{L}\right),$$

$$\bar{w}(T,x) = T^{-1} \left(-\frac{8\pi}{3L} \cos\left(\frac{\pi x}{L}\right) \sin\left(\frac{\pi x}{L}\right)\right),$$

$$\forall k \in \{0, \dots, N-1\}, \ \forall x \in \mathcal{C}_k,$$

$$\bar{u}(t_k,x) = \mathbb{E}[\bar{u}(t_{k+1}, \Pi_{k+1}(x + \bar{w}(t_{k+1},x)h + \sqrt{2\bar{u}(t_{k+1},x)} g(\Delta B^k)))]$$
$$\qquad\qquad + h\bar{u}(t_{k+1},x)^2,$$

$$\bar{w}(t_k,x) = \mathbb{E}[\bar{w}(t_{k+1}, \Pi_{k+1}(x + 3\bar{w}(t_{k+1},x)h + \sqrt{2\bar{u}(t_{k+1},x)} g(\Delta B^k)))]$$
$$\qquad\qquad + 2h\bar{u}(t_{k+1},x)\bar{w}(t_{k+1},x).$$

For $T = 1, h = 0.02, \ \delta = L/500, \ M = 160$, we present below the results obtained first with Algorithm 3.1 (the approximation of the gradient with this algorithm is undefined at $x = \pm L/2$ and we thus arbitrarily set it to zero) and then with Algorithm 5.3. See Figure 5 for the results on $[-L/2, L/2]$.

We first observe that the approximated solutions obtained with the two algorithms are not significantly different. The main advantage of the differentiated algorithm is, as expected, for the pointwise approximation of the gradient. Indeed, in that case there is a factor 4 between the absolute pointwise errors associated to the two methods. Let us also indicate that both methods present some "singularity" in the neighborhood of $x = \pm L/2$ for the estimation of the gradient. This could be expected for Algorithm 3.1 since the estimation of the gradient is obtained by dividing $\bar{v}$ by $\sqrt{2\bar{u}}$ that goes to 0 when $x \to \pm L/2$. It is a bit more surprising for Algorithm 5.3.

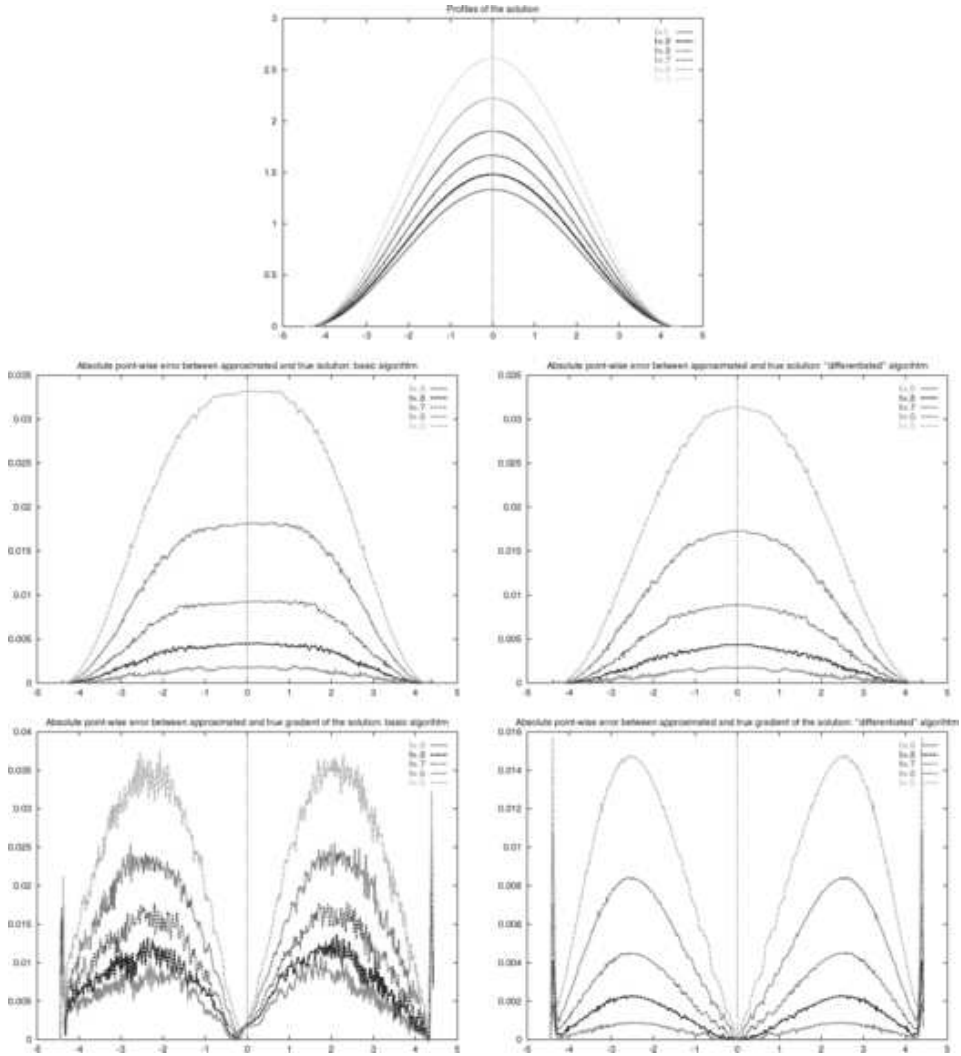

Fig. 5.

**6. Proof. First step: a priori controls.** In this section we give various a priori estimates of the couple $(Y, Z)$ introduced in (3.13) and of the approximate diffusion $X$ defined in (3.12). These controls are necessary to establish Theorems 4.1, 4.2 and 4.3.

*About constants.* In the following, we keep the same notation $C, C_\vartheta, c_\vartheta$ (or $C', C'_\vartheta, c'_\vartheta$) for all finite, nonnegative constants which appear in our computations: they may depend on known parameters deriving from Assumption (A), on $T$ and on $p$, but not on any of the discretization parameters. The index $\vartheta$ in the previous notation refers to the numbering of the Proposition, Lemma, Theorem, . . . where the constant appears.

*Conditions on parameters.* We assume that the conditions of Theorem 4.1 on $h$, $\delta$, $M$, $\rho$ and $p$ are fulfilled.

6.1. *Discrete backward equation and a priori estimates.*

*Discrete Feynman–Kac formula.* By iteration of the dynamic programming principle in Algorithm 3.1, it is plain to prove the discrete Feynman–Kac formula (3.14).

Both formulae (3.14) and (3.15) [representation of $Y_{t_N} + h \sum_{i=1}^{N} f(X_{t_{i-1}}, \bar{u}(t_i, X_{t_{i-1}}), Z_{t_{i-1}})$ through the martingale representation theorem] permit to apply the BSDE machinery to our frame. However, as well known in the literature devoted to SDEs (or, equivalently, to PDEs), several a priori estimates of the solution are necessary to apply this strategy.

PROPOSITION 6.1. *There exists a constant $C_{6.1}$ s.t.*

$$\sup_{i=0,\ldots,N} \left[ \sup_{x \in \mathcal{C}_i} |\bar{u}(t_i, x)|^2 \right] \leq C_{6.1}.$$

PROPOSITION 6.2. *There exists a constant $C_{6.2}$ s.t.*

$$\mathbb{E}\left[ \int_0^T |\overline{Z}_s|^2 \, ds \right] + h \sum_{i=0}^{N-1} \mathbb{E}[|Z_{t_i}|^2]$$

$$+ h \sup_{i=0,\ldots,N} \left[ \sup_{x \in \mathcal{C}_i} |\bar{v}(t_i, x)|^2 \right] + h \sup_{i=0,\ldots,N-1} \left[ \sup_{x \in \mathcal{C}_i} |\hat{v}(t_i, x)|^2 \right] \leq C_{6.2}.$$

The distance between $Z$ and $\overline{Z}$ can be estimated as follows:

LEMMA 6.3. *There exists a constant $C_{6.3}$ s.t., for $k \in \{1, \ldots, N\}$,*

$$\mathbb{E}\left| hZ_{t_{k-1}} - \mathbb{E}\left[ \int_{t_{k-1}}^{t_k} \overline{Z}_s \, ds \Big| \mathcal{F}_{t_{k-1}} \right] \right|^2 \leq C_{6.3} h^2 \mathcal{E}^2(\text{quantiz}).$$

### 6.2. Approximate diffusion.

*Jumps of the discrete forward process.*   Start first with the following:

LEMMA 6.4.   *For a given $k \in \{0, \ldots, N-1\}$, the norm of the increment $X_{t_{k+1}} - X_{t_k}$ is always bounded by $|\mathcal{T}(t_k, X_{t_k})| + \delta$. In particular, there exists a constant $C_{6.4}$ such that*

$$\mathbb{E}[|X_{t_{k+1}} - X_{t_k}|^2 | \mathcal{F}_{t_k}] \leq C_{6.4}[h + \delta^2].$$

PROOF.   Since $X_{t_k} \in \mathcal{C}_\infty$, one has $\Pi_\infty(X_{t_k} + \mathcal{T}(t_k, X_{t_k})) = X_{t_k} + \Pi_\infty(\mathcal{T}(t_k, X_{t_k}))$ (invariance by translation of the grid $\mathcal{C}_\infty$). Moreover, for every $y$ in the image of the projection $\mathcal{Q}(R + \rho, \cdot)$ and for every $z \in \mathbb{R}^d$, the distance $|\mathcal{Q}(R + \rho, y + z) - y|$ is bounded by $|z|$. Hence,

$$(6.1) \quad \begin{aligned} |X_{t_{k+1}} - X_{t_k}| &= |\mathcal{Q}(R + \rho, X_{t_k} + \Pi_\infty(\mathcal{T}(t_k, X_{t_k}))) - X_{t_k}| \\ &\leq |\Pi_\infty(\mathcal{T}(t_k, X_{t_k}))| \leq |\mathcal{T}(t_k, X_{t_k})| + \delta. \end{aligned}$$

Thanks to Propositions 6.1 and 6.2, we are able to bound the drift $b$ appearing in the transition. Since $\mathbb{E}[|g(\Delta B^k)|^2] \leq Ch$, from Assumption (A) and Proposition 6.1, we also control the martingale part of the transition. This completes the proof.   □

*Extension of the "discrete diffusion".*   For the proof, we need to extend the definition of $X$ to the whole set $[0, T]$. Put, for all $k \in \{0, \ldots, N-1\}$ and $t \in [t_k, t_{k+1}[$,

$$(6.2) \quad \begin{aligned} X_t &\equiv X_{t_k} + b(X_{t_k}, \bar{u}(t_{k+1}, X_{t_k}), \hat{v}(t_k, X_{t_k}))(t - t_k) \\ &\quad + \sigma(X_{t_k}, \bar{u}(t_{k+1}, X_{t_k}))[B_t - B_{t_k}]. \end{aligned}$$

From Proposition 6.2, we get the following:

LEMMA 6.5.   *There exists a constant $C_{6.5}$ s.t., for every $k \in \{0, \ldots, N-1\}$,*

$$\forall t \in [t_k, t_{k+1}[, \qquad \mathbb{E}[|X_t - X_{t_k}|^2 | \mathcal{F}_{t_k}] \leq C_{6.5}h.$$

The extended process $(X_t)_{0 \leq t \leq T}$ is discontinuous at times $(t_k)_{1 \leq k \leq N}$. At a given time $t_k$, $1 \leq k \leq N$, the size of the jump performed by the process depends on the quantization error and on the spatial projection error. The first error is easily controlled by the distortion. Concerning the second one, the projection error is close to the spatial step $\delta$ when the grids are infinite. For truncated grids, the story is slightly different. In fact, as soon as the process stays inside $(\Delta_k)_{0 \leq k \leq N}$, the projection error is close to the step $\delta$ of the interior mesh of the grid $(\mathcal{C}_k)_{0 \le k \le N}$. At the opposite, outside $(\Delta_k)_{0 \le k \le N}$, the jump of the process may take large values.

The time continuous extension of $X$ remains close to the discrete version of $X$ up to time $\tau_\infty$.

LEMMA 6.6. *There exists a constant* $C_{6.6}$ *such that*

$$\sum_{i=0}^{N-1} \mathbb{E}[\mathbf{1}_{\{t_{i+1} < \tau_\infty\}}|X_{t_{i+1}} - X_{t_{i+1}-}|^2] \le C_{6.6}h(\mathcal{E}^2(\text{space}) + \mathcal{E}^2(\text{quantiz})).$$

PROOF (SKETCH). From (6.2), the difference $X_{t_{i+1}} - X_{t_{i+1}-}$ writes

$$
\begin{aligned}
X_{t_{i+1}} - X_{t_{i+1}-} = {}& [\Pi_{i+1}(X_{t_i} + \mathcal{T}(t_i, X_{t_i})) - (X_{t_i} + \mathcal{T}(t_i, X_{t_i}))] \\
& + \sigma(X_{t_i}, \bar{u}(t_{i+1}, X_{t_i}))[g(\Delta B^i) - \Delta B^i] \\
\equiv {}& E_1(i+1) + E_2(i+1).
\end{aligned}
$$
(6.3)

$E_1(i+1)$ appears as a projection error and $E_2(i+1)$ as a quantization one. It is readily seen that $E_1(i+1)$ is bounded by $\delta$ on $\{t_{i+1} < \tau_\infty\}$. From (3.8), one also gets $\mathbb{E}[|E_2(i+1)|^2|\mathcal{F}_{t_i}] \le ChM^{-2/d}$. □

### 6.3. *Sketches of the proofs of the a priori controls.*

*Discrete BSDE.* This section is devoted to the proof of Propositions 6.1, 6.2 and Lemma 6.3. We first give a control of the $L^2$ norm between $Z_{t_{k-1}}$ and the conditional expectation of $\int_{t_{k-1}}^{t_k} \overline{Z}_s \, ds$ appearing in Lemma 6.3. This preliminary estimate permits to prove Proposition 6.1. We then derive the complete proofs of Proposition 6.2 and Lemma 6.3.

*Step one*: *preliminary control in Lemma* 6.3. From (3.15), write, for a given $k \in \{0, \ldots, N-1\}$,

$$Y_{t_{k+1}} + hf(X_{t_k}, \bar{u}(t_{k+1}, X_{t_k}), Z_{t_k}) = Y_{t_k} + \int_{t_k}^{t_{k+1}} \overline{Z}_s \, dB_s.$$

Multiply this identity by $\Delta B^k$, take the conditional expectation w.r.t. $\mathcal{F}_{t_k}$ and plug the definition of $Z_{t_k}$ [cf. (3.13)]:

$$hZ_{t_k} - \mathbb{E}\left[\int_{t_k}^{t_{k+1}} \overline{Z}_s \, ds \Big| \mathcal{F}_{t_k}\right] = \mathbb{E}[Y_{t_{k+1}}(g(\Delta B^k) - \Delta B^k)|\mathcal{F}_{t_k}].$$
(6.4)

Referring to (3.8), there exists $C$ s.t.

$$\mathbb{E}\left[\left|hZ_{t_k} - \mathbb{E}\left[\int_{t_k}^{t_{k+1}} \overline{Z}_s \, ds \Big| \mathcal{F}_{t_k}\right]\right|^2\right] \le ChM^{-2/d}\mathbb{E}[Y_{t_{k+1}}^2].$$
(6.5)

This preliminary estimate (6.5) is necessary to prove Proposition 6.1 from which we will derive $\mathbb{E}[Y_{t_{k+1}}^2] \le C$, and thus complete the proof of Lemma 6.3.

*Step two*: *proof of Proposition* 6.1.   To estimate the *supremum* norm of $\bar{u}$ over the grids $\mathcal{C}_0, \ldots, \mathcal{C}_N$, we follow the basic strategy of the BSDE theory and, therefore, apply a discrete version of Itô's formula to the discrete BSDE formula given in (3.14)–(3.15). Such a formula can be found in [22], Chapter VII, Section 9. We obtain

$$
\mathbb{E}|Y_T|^2 = |Y_0|^2 + 2h \sum_{i=1}^{N} \mathbb{E}\langle -f(X_{t_{i-1}}, \bar{u}(t_i, X_{t_{i-1}}), Z_{t_{i-1}}), Y_{t_{i-1}} \rangle
$$
$$
(6.6)
$$
$$
+ h^2 \sum_{i=1}^{N} \mathbb{E}[f^2(X_{t_{i-1}}, \bar{u}(t_i, X_{t_{i-1}}), Z_{t_{i-1}})] + \mathbb{E} \int_0^T |\overline{Z}_s|^2 \, ds.
$$

Following standard computations in BSDE theory, it is plain to derive from (6.5) and (6.6):

$$
(6.7) \quad |\bar{u}(0, x_0)|^2 + \mathbb{E} \int_0^T |\overline{Z}_s|^2 \, ds + h \sum_{i=0}^{N-1} \mathbb{E}[|Z_{t_i}|^2] \le C + Ch \sum_{i=0}^{N} \sup_{x \in \mathcal{C}_i} |\bar{u}(t_i, x)|^2.
$$

There exists a constant $c > 0$ such that, for $h < c$ (recall indeed that $h$ is small), the above inequality holds but with $i = 1$ instead of $i = 0$ as initial condition in the r.h.s. of (6.7). As usual in BSDE theory, we can establish in a similar way that, for every initial condition $(t_k, x)$, $1 \le k \le N$,

$$
\forall k \in \{0, \ldots, N-1\}, \qquad \sup_{x \in \mathcal{C}_k} |\bar{u}(t_k, x)|^2 \le C + Ch \sum_{i=k+1}^{N} \sup_{x \in \mathcal{C}_i} |\bar{u}(t_i, x)|^2.
$$

A discrete version of Gronwall's lemma yields the result.

*Step three*: *proofs of Proposition* 6.2 *and Lemma* 6.3.   The $L^2$-estimates of $Z$ and $\overline{Z}$ in Proposition 6.2 follow from Proposition 6.1 and (6.7). Moreover, as a consequence of Proposition 6.1 and the definitions of $\bar{v}$ and $\hat{v}$, see Algorithm 3.1, we deduce the estimates of the *supremum* norms of $\bar{v}$ and $\hat{v}$. Lemma 6.3 follows from (6.5) and Proposition 6.1.

## 7. Proof. Second step: stability properties.   This section focuses on the second step of the proof of Theorems 4.1, 4.2 and 4.3, and aims to establish more specifically a suitable intermediate inequality, close to usual stability properties of FBSDEs.

*Strategy*.   Recall first that two main strategies are conceivable in the theoretical framework to establish classical stability theorems for FBSDEs.

Denote to this end by $(U', V', W')$ a solution of another FBSDE of type $(E)$ with different coefficients. The associated PDE solution is just denoted by $u'$. In order to compare $u'$ with $u$, the following approaches have been employed in the literature:

1. First, the recent induction principle given in [6] can be applied. In short, $u$ and $u'$ are compared on a neighborhood of the boundary $T$ with classical arguments of stochastic analysis and the estimate of the difference between these solutions is then extended by induction from the final bound $T$ to the initial bound 0. The local estimates consist in studying the distance between $U$ and $U'$ and between $(V, W)$ and $(V', W')$. This strategy has been successfully applied to various contexts (see [6] for the solvability of FBSDEs and [8] for homogenization of quasilinear PDEs).

2. A second approach follows the earlier *Four Step Scheme* of Ma, Protter and Yong [14]. Instead of studying the difference between $U$ and $U'$ and between $(V, W)$ and $(V', W')$, the process $(u(t, U'_t))_{0 \leq t \leq T}$ is written with Itô's formula as the solution of a BSDE. This BSDE is then compared with the one satisfied by $(V', W')$. In particular, these BSDEs are both written with respect to the same diffusion $U'$. Generally speaking, this strategy holds when $u$ is smooth enough (e.g., if $u$ satisfies Theorem 2.1). It is then more direct than the previous one.

Under Assumption (A) we apply the second strategy and compare the process $Y$ with the process $(u(t, X_t))_{0 \leq t \leq T \wedge \tau_\infty}$ [see (6.2) for the definition of the extension of $X$].

7.1. *Statements of the stability results.*

*First stability property.* Applying the usual FBSDE machinery, we are able to establish in Section 7.2 the following first inequality:

PROPOSITION 7.1. *There exists a constant $C_{7.1}$ such that, for $\eta$ small enough,*

$$|(\bar{u} - u)(0, x_0)|^2 + C_{7.1}^{-1} h \sum_{j=1}^{N} \mathbb{E}[|(\bar{v} - v)(t_{j-1}, X_{t_{j-1}})|^2 \mathbf{1}_{\{t_{j-1} < \tau_\infty\}}]$$

$$\leq C_{7.1} \Bigg[ \mathbb{P}\{\tau_\infty < +\infty\} + \mathcal{E}^2(\text{time}) + \mathcal{E}^2(\text{space}) + \mathcal{E}^2(\text{quantiz})$$

$$(7.1) \qquad + \eta^{-1} h \sum_{j=1}^{N} \mathbb{E}[|(\bar{u} - u)(t_j, X_{t_{j-1}})|^2 \mathbf{1}_{\{t_{j-1} < \tau_\infty\}}]$$

$$+ \eta^{-1} h \sum_{j=1}^{N} \mathbb{E}[|(\bar{u} - u)(t_{j-1}, X_{t_{j-1}})|^2 \mathbf{1}_{\{t_{j-1} < \tau_\infty\}}]$$

$$+ (\eta + h) h \sum_{j=1}^{N} \mathbb{E}[|(\hat{v} - v)(t_{j-1}, X_{t_{j-1}})|^2 \mathbf{1}_{\{t_{j-1} < \tau_\infty\}}] \Bigg].$$

*When the drift b does not depend on z, the last term of the r.h.s. does not appear.*

*Estimates of the gradient increment.* Assume for the moment that Proposition 7.1 holds. Note that the main problem then remains to estimate the last term in the r.h.s. of (7.1). Thanks to the specific choice of $\hat{v}$ in Section 3.3, we are able to establish in Section 7.3 the following control:

PROPOSITION 7.2. *There exists a constant $C_{7.2}$ such that, for $k \in \{0, \dots, N-1\}$, on $\{t_k < \tau_\infty\}$,*

$$|(\hat{v} - v)(t_k, X_{t_k})| \leq C_{7.2}[\mathcal{E}(\text{gradient}, p) + \mathcal{E}(\text{time}) + h\mathcal{E}(\text{space})$$
$$+ \mathbb{E}[|(\bar{v} - v)(t_{k+1}, X_{t_{k+1}})|^2 | \mathcal{F}_{t_k}]^{1/2}].$$

*Main stability theorem.* From Propositions 7.1 and 7.2, we claim the following:

THEOREM 7.3. *Proposition 7.1 holds with the last term in the r.h.s. of (7.1) replaced by*

$$\mathcal{E}^2(\text{gradient}, p) + (\eta + h)h \sum_{j=1}^{N} \mathbb{E}[|(\bar{v} - v)(t_j, X_{t_j})|^2 \mathbf{1}_{\{t_{j-1} < \tau_\infty\}}].$$

Application of Theorem 7.3 to the proof of Theorems 4.1, 4.2 and 4.3 is given in Section 8.

### 7.2. *Proof of Proposition 7.1.*

*Starting point*: *time continuous backward processes.* Following the second strategy and referring to the structure of the PDE ($\mathcal{E}$), set for notational convenience

$$(7.2) \quad \forall t \in [0, T], \qquad \overline{V}_t \equiv u(t, X_t), \qquad \overline{W}_t \equiv \nabla_x u(t, X_t)\sigma(X_t, \overline{V}_t).$$

Note, moreover, that the martingale part of $(\overline{V}_t)_{0 \leq t \leq T}$ is driven by

$$(7.3) \quad \forall t \in [0, T[, \qquad \hat{W}_t \equiv \nabla_x u(t, X_t)\sigma(X_{\phi(t)}, \bar{u}(\phi(t) + h, X_{\phi(t)})),$$

where $\phi(t) = t_k$ for $t_k \leq t < t_{k+1}$, $k \in \{0, \dots, N-1\}$. From Theorem 2.1 and Lemma 6.5, we derive the following *a priori estimates* of $\overline{V}, \overline{W}$ for $s \in [t_k, t_{k+1}[$:

$$(7.4) \qquad \mathbb{E}[|\overline{V}_s - \overline{V}_{t_k}| + |\overline{W}_s - \overline{W}_{t_k}| | \mathcal{F}_{t_k}] \leq Ch^{1/2}.$$

*Step one*: *Itô's formula for* $\overline{V}$.   Using Itô's formula and the PDE satisfied by $u$, we obtain, for $i \in \{0, \ldots, N-1\}$,

$$\overline{V}_{t_{i+1}} - \overline{V}_{t_i} = \overline{V}_{t_{i+1}} - \overline{V}_{t_{i+1}-} + \int_{t_i}^{t_{i+1}} [F(s, X_s, X_{t_i}, \bar{u}(t_{i+1}, X_{t_i}), \hat{v}(t_i, X_{t_i}))$$
$$- F(s, X_s, X_s, \overline{V}_s, \overline{W}_s)] \, ds$$
$$- \int_{t_i}^{t_{i+1}} f(X_s, \overline{V}_s, \overline{W}_s) \, ds + \int_{t_i}^{t_{i+1}} \hat{W}_s \, dB_s,$$

with $F(s, x, \hat{x}, y, z) = \langle \nabla_x u(s, x), b(\hat{x}, y, z) \rangle + (1/2) \operatorname{tr}(a(\hat{x}, y) \nabla_{x,x}^2 u(s, x))$.

*Step two*: *difference of the processes*.   The strategy is well known: we aim to make the difference between $\overline{V}$ and $Y$ and then to apply the usual BSDE machinery to estimate the distance between these processes. Hence, we claim from (3.15)

$$\overline{V}_{t_{i+1}} - Y_{t_{i+1}} - [\overline{V}_{t_i} - Y_{t_i}]$$
$$= \overline{V}_{t_{i+1}} - \overline{V}_{t_{i+1}-}$$
$$+ \int_{t_i}^{t_{i+1}} [F(s, X_s, X_{t_i}, \bar{u}(t_{i+1}, X_{t_i}), \hat{v}(t_i, X_{t_i})) - F(s, X_s, X_s, \overline{V}_s, \overline{W}_s)] \, ds$$
$$- \int_{t_i}^{t_{i+1}} [f(X_s, \overline{V}_s, \overline{W}_s) - f(X_{t_i}, \bar{u}(t_{i+1}, X_{t_i}), Z_{t_i})] \, ds$$
$$+ \int_{t_i}^{t_{i+1}} [\hat{W}_s - \overline{Z}_s] \, dB_s$$
$$\equiv \Delta E_{i+1}(1) + \Delta E_{i+1}(2) + \Delta E_{i+1}(3) + \Delta E_{i+1}(4).$$

The discrete Itô formula [see the derivation of (6.6)] and standard computations yield

$$(7.5) \qquad |\overline{V}_0 - Y_0|^2 + \tfrac{1}{2} D(3) \leq \mathbb{E} |\overline{V}_{T \wedge \tau_\infty} - Y_{T \wedge \tau_\infty}|^2 + D(1) + D(2),$$

with

$$D(1) \equiv -2 \mathbb{E} \sum_{j=1}^{N} [\mathbf{1}_{\{t_{j-1} < \tau_\infty\}} [\overline{V}_{t_{j-1}} - Y_{t_{j-1}}] E_j],$$

$$(7.6) \quad D(2) \equiv \sum_{j=1}^{N} \mathbb{E}[\mathbf{1}_{\{t_{j-1} < \tau_\infty\}} E_j^2], \qquad D(3) \equiv \sum_{j=1}^{N} \mathbb{E}[\mathbf{1}_{\{t_{j-1} < \tau_\infty\}} \Delta E_j(4)^2],$$

$$E_j \equiv \Delta E_j(1) + \Delta E_j(2) + \Delta E_j(3), \qquad j \in \{1, \ldots, N\}.$$

*Step three*: *standard BSDE techniques*.   Following the BSDE techniques, we have to upper bound $D(1), D(2)$ [resp. lower bound $D(3)$] by terms ap-pearing in the r.h.s. (resp. l.h.s.) of (7.1). The following lemmas whose proofs are postponed to the end of the subsection give the needed controls.

LEMMA 7.4.   *Denote by* RHS(7.1) *the r.h.s. of* (7.1). *Then, there exists a constant* $C_{7.4}$ *such that, for* $\eta \in \, ]0,1]$,

$$|D(1)| + D(2) \leq C \left[ \mathrm{RHS}(7.1) \right.$$
$$\left. + h(\eta + h) \sum_{j=1}^{N} \mathbb{E}[|(\bar{v} - v)(t_{j-1}, X_{t_{j-1}})|^2 \mathbf{1}_{\{t_{j-1} < \tau_\infty\}}] \right].$$

LEMMA 7.5.   *There exists a constant* $C_{7.5} > 0$ *such that*

$$D(3) \geq C_{7.5}^{-1} h \sum_{j=1}^{N} \mathbb{E}[\mathbf{1}_{\{t_{j-1} < \tau_\infty\}} |(\bar{v} - v)(t_{j-1}, X_{t_{j-1}})|^2]$$
$$- C_{7.5}(\mathcal{E}^2(\text{quantiz}) + \mathcal{E}^2(\text{time}))$$
$$- C_{7.5} h \sum_{j=1}^{N} \mathbb{E}[\mathbf{1}_{\{t_{j-1} < \tau_\infty\}} |(\bar{u} - u)(t_j, X_{t_{j-1}})|^2].$$

Note to conclude the proof of Proposition 7.1 that $Y_T = \overline{V}_T$. Hence, from Theorem 2.1 and Proposition 6.1 (boundedness of $u$ and $\bar{u}$), $\mathbb{E}|\overline{V}_{T \wedge \tau_\infty} - Y_{T \wedge \tau_\infty}|^2 \leq C\mathbb{P}\{\tau_\infty < T\} \leq C\mathbb{P}\{\tau_\infty < +\infty\}$. Choose finally $\eta$ small enough to obtain inequality (7.1) from (7.5), (7.6), and Lemmas 7.4 and 7.5. This completes, up to the proofs of Lemmas 7.4 and 7.5, the proof of Proposition 7.1.

PROOF OF LEMMA 7.4.   Note from Theorem 2.1 that $\Delta E_j(2)$ and $\Delta E_j(3)$ may be seen as "Lipschitz" differences since the partial derivatives of $u$ of order one and two in $x$ are bounded. Recall also that $\overline{V}_s = u(s, X_s)$, $\overline{W}_{t_{j-1}} = v(t_{j-1}, X_{t_{j-1}})$ and $Z_{t_{j-1}} = \bar{v}(t_{j-1}, X_{t_{j-1}})$. From Theorem 2.1 (Hölder regularity of $u$ in $t$), (7.4) (regularity of $\overline{V}$ and $\overline{W}$), Lemma 6.5 (control of the increments of $X$) and Young's inequality, it comes, for every $\eta \in \, ]0,1]$,

$$|D(1)| \leq C\mathcal{E}^2(\text{time})$$
$$+ C\mathbb{E} \sum_{j=1}^{N} [\mathbf{1}_{\{t_{j-1} < \tau_\infty\}} |\overline{V}_{t_{j-1}} - Y_{t_{j-1}}| |\overline{V}_{t_j} - \overline{V}_{t_{j-1}}|]$$
$$(7.7) \qquad + Ch \sum_{j=1}^{N} [\eta^{-1} \mathbb{E}[|(\bar{u} - u)(t_{j-1}, X_{t_{j-1}})|^2 \mathbf{1}_{\{t_{j-1} < \tau_\infty\}}]$$

$$+ \mathbb{E}[|(\bar{u} - u)(t_j, X_{t_{j-1}})|^2 \mathbf{1}_{\{t_{j-1} < \tau_\infty\}}]]$$

$$+ \eta h \sum_{j=1}^{N} [\mathbb{E}[|(\hat{v} - v)(t_{j-1}, X_{t_{j-1}})|^2 \mathbf{1}_{\{t_{j-1} < \tau_\infty\}}]$$

$$+ \mathbb{E}[|(\bar{v} - v)(t_{j-1}, X_{t_{j-1}})|^2 \mathbf{1}_{\{t_{j-1} < \tau_\infty\}}]].$$

It now remains to estimate the second term in the r.h.s. of (7.7). Note first that, for all $j \in \{1, \ldots, N\}$, $\{t_{j-1} < \tau_\infty\} = \{t_j < \tau_\infty\} \cup \{t_j = \tau_\infty\}$. Hence, thanks to the boundedness of $u$ and $\bar{u}$ (see Theorem 2.1 and Proposition 6.1), to Lemma 6.6 (jumps of the process $X$) and to the global Lipschitz property of $u$ (see Theorem 2.1),

$$\mathbb{E} \sum_{j=1}^{N} [\mathbf{1}_{\{t_{j-1} < \tau_\infty\}} [|\overline{V}_{t_{j-1}} - Y_{t_{j-1}}||\overline{V}_{t_j} - \overline{V}_{t_j-}|]]$$

$$\leq \mathbb{E} \sum_{j=1}^{N} [\mathbf{1}_{\{t_j < \tau_\infty\}} [|\overline{V}_{t_{j-1}} - Y_{t_{j-1}}|^2 + |\overline{V}_{t_j} - \overline{V}_{t_j-}|^2]]$$

(7.8) $$+ C\mathbb{P}\{\tau_\infty < +\infty\}$$

$$\leq Ch\mathbb{E} \sum_{j=1}^{N} [\mathbf{1}_{\{t_j < \tau_\infty\}} |(\bar{u} - u)(t_{j-1}, X_{t_{j-1}})|^2]$$

$$+ C(\mathcal{E}^2(\text{space}) + \mathcal{E}^2(\text{quantiz})) + C\mathbb{P}\{\tau_\infty < +\infty\}.$$

Plug (7.8) in (7.7) to derive the required control for $D(1)$.

Turn to the estimation of $D(2)$: apply again Lemma 6.6 to control $\Delta E_j(1)$, for a given $j \in \{1, \ldots, N\}$, and treat the "Lipschitz" differences as done to estimate $D(1)$. □

PROOF OF LEMMA 7.5.    Write first

$$h \sum_{j=1}^{N} \mathbb{E}[\mathbf{1}_{\{t_{j-1} < \tau_\infty\}} |\bar{v}(t_{j-1}, X_{t_{j-1}}) - v(t_{j-1}, X_{t_{j-1}})|^2]$$

$$\leq Ch \sum_{j=1}^{N} \left\{ \mathbb{E}\left[ \mathbf{1}_{\{t_{j-1} < \tau_\infty\}} \left| \bar{v}(t_{j-1}, X_{t_{j-1}}) - \frac{1}{h}\mathbb{E}\left[ \int_{t_{j-1}}^{t_j} \overline{Z}_s \, ds \Big| \mathcal{F}_{t_{j-1}} \right] \right|^2 \right] \right.$$

$$+ \mathbb{E}\left[ \mathbf{1}_{\{t_{j-1} < \tau_\infty\}} \left| \frac{1}{h}\mathbb{E}\left[ \int_{t_{j-1}}^{t_j} [\overline{Z}_s - \hat{W}_s] \, ds \Big| \mathcal{F}_{t_{j-1}} \right] \right|^2 \right]$$

$$\left. + \mathbb{E}\left[ \mathbf{1}_{\{t_{j-1} < \tau_\infty\}} \left| \frac{1}{h}\mathbb{E}\left[ \int_{t_{j-1}}^{t_j} [\hat{W}_s - v(t_{j-1}, X_{t_{j-1}})] \, ds \Big| \mathcal{F}_{t_{j-1}} \right] \right|^2 \right] \right\}$$

$$\equiv A(1) + A(2) + A(3).$$

From Lemma 6.3 (distance between $Z$ and $\overline{Z}$), we then derive $A(1) \leq C\mathcal{E}^2(\text{quantiz})$. For the term $A(2)$, the Cauchy–Schwarz inequality yields $A(2) \leq CD(3)$. Concerning $A(3)$, we get

$$
A(3) \leq C \sum_{j=1}^{N} \mathbb{E}\bigg[ \mathbf{1}_{\{t_{j-1} < \tau_\infty\}}
$$
$$
\times \int_{t_{j-1}}^{t_j} |\nabla_x u(s, X_s)\sigma(X_{t_{j-1}}, \bar{u}(t_j, X_{t_{j-1}}))
$$
$$
- \nabla_x u(t_{j-1}, X_{t_{j-1}})\sigma(X_{t_{j-1}}, u(t_{j-1}, X_{t_{j-1}}))|^2 \, ds \bigg].
$$

Following the techniques employed in the previous proof, relying on the smoothness of the true solution (see Theorem 2.1) on the boundedness of the approximate solution, see Proposition 6.1, and on intermediate controls of the process $X$, see Lemma 6.5, we get

$$
A(3) \leq Ch\bigg[ 1 + \sum_{j=1}^{N} \mathbb{E}[\mathbf{1}_{\{t_{j-1} < \tau_\infty\}} |\bar{u}(t_j, X_{t_{j-1}}) - u(t_j, X_{t_{j-1}})|^2] \bigg].
$$

The above estimates of $A(1), A(2), A(3)$ complete the proof. $\quad\square$

### 7.3. *Proof of Proposition* 7.2 (*difference of the gradients*).

*Strategy.* In Proposition 7.2, we aim to control the quantity $|(\hat{v} - v)(t_k, X_{t_k})|$ for $t_k < \tau_\infty$, with $\hat{v}(t_k, X_{t_k}) = \mathbb{E}[\bar{v}(t_{k+1}, \Pi_{k+1}(X_{t_k} + \mathcal{T}^0(t_k, X_{t_k})))|\mathcal{F}_{t_k}]$ (see Algorithm 3.1). We first write $v(t_k, X_{t_k})$ in a similar way to study the difference $(\hat{v} - v)(t_k, X_{t_k})$. From Theorem 2.1 (regularity of $u$) and from the proof of Lemma 6.4 (with $\mathcal{T}^0$ instead of $\mathcal{T}$), we claim

$$
|\mathbb{E}[v(t_{k+1}, \Pi_{k+1}(X_{t_k} + \mathcal{T}^0(t_k, X_{t_k})))|\mathcal{F}_{t_k}] - v(t_k, X_{t_k})| \leq C[h^{1/2} + \delta].
$$

Hence,

$$
\begin{aligned}
(7.9) \quad &|(\hat{v} - v)(t_k, X_{t_k})| \\
&\leq C\mathbb{E}[|(\bar{v} - v)(t_{k+1}, \Pi_{k+1}(X_{t_k} + \mathcal{T}^0(t_k, X_{t_k})))||\mathcal{F}_{t_k}] + C(h^{1/2} + \delta).
\end{aligned}
$$

Proposition 7.2 directly follows from (7.9) and the next theorem:

THEOREM 7.6. *There exists a constant* $C_{7.6}$ *such that on* $\{t_k < \tau_\infty\}$

$$
\mathbb{E}[|(\bar{v} - v)(t_{k+1}, \Pi_{k+1}(X_{t_k} + \mathcal{T}^0(t_k, X_{t_k})))||\mathcal{F}_{t_k}]
$$
$$
\leq C_{7.6}\mathcal{E}(\text{gradient}, p) + C_{7.6}\mathbb{E}[|(\bar{v} - v)(t_{k+1}, X_{t_{k+1}})|^2|\mathcal{F}_{t_k}]^{1/2}.
$$

The main difficulty to prove Theorem 7.6 lies in the lack of regularity of $\bar{v}$. To overcome this point, note first that

$$(7.10) \qquad \mathbb{E}[|(\bar{v} - v)(t_{k+1}, \Pi_{k+1}(X_{t_k} + \mathcal{T}^0(t_k, X_{t_k})))||\mathcal{F}_{t_k}]$$

and

$$(7.11) \qquad \mathbb{E}[|(\bar{v} - v)(t_{k+1}, \Pi_{k+1}(X_{t_k} + \mathcal{T}(t_k, X_{t_k})))|^2|\mathcal{F}_{t_k}]^{1/2}$$

write as expectations of a given function with respect to two different kernels. We then aim to compare these underlying kernels. Recall that for a given $x \in \mathcal{C}_k$, both $\mathcal{T}^0(t_k, x)$ and $\mathcal{T}(t_k, x)$ are, up to a quantization procedure, Gaussian random variables with same covariance matrices but different means. The strategy then consists in applying a Gaussian change of variable to pass from the first kernel to the second one.

*Step one*: *Proof of Theorem 7.6, exhibition of underlying kernels.* We first write (7.10) with respect to the underlying kernel $\mathcal{T}^0$. Note in this frame, with the notation of Section 3.4, that, for every $x \in \mathbb{R}^d$, $\Pi_{k+1}(x) = \Pi_{k+1} \circ \Pi_\infty(x)$ since $\Pi_\infty(x) \in \Delta_k \Leftrightarrow x \in \Delta_k$. Thus, using the invariance by translation of $\mathcal{C}_\infty$ (see the proof of Lemma 6.4), (7.10) writes

$$(7.12) \qquad \begin{aligned} &\mathbb{E}[|(\bar{v} - v)(t_{k+1}, \Pi_{k+1}(X_{t_k} + \mathcal{T}^0(t_k, X_{t_k})))||\mathcal{F}_{t_k}] \\ &= \sum_{y \in \mathcal{C}_\infty} [|(\bar{v} - v)(t_{k+1}, \Pi_{k+1}(X_{t_k} + y))| \\ &\qquad \times \mathbb{P}\{\Pi_\infty(\mathcal{T}^0(t_k, X_{t_k})) = y|\mathcal{F}_{t_k}\}]. \end{aligned}$$

In the same way, the square of (7.11) writes

$$(7.13) \qquad \begin{aligned} &\mathbb{E}[|(\bar{v} - v)(t_{k+1}, \Pi_{k+1}(X_{t_k} + \mathcal{T}(t_k, X_{t_k})))|^2|\mathcal{F}_{t_k}] \\ &= \sum_{y \in \mathcal{C}_\infty} [|(\bar{v} - v)(t_{k+1}, \Pi_{k+1}(X_{t_k} + y))|^2 \\ &\qquad \times \mathbb{P}\{\Pi_\infty(\mathcal{T}(t_k, X_{t_k})) = y|\mathcal{F}_{t_k}\}]. \end{aligned}$$

Equations (7.12) and (7.13) provide relevant writings to estimate (7.10) and (7.11). Indeed, it is sufficient to bound for a given $x \in \mathcal{C}_k$ and a given $y \in \mathcal{C}_\infty$ the probability $\mathbb{P}\{\Pi_\infty(\mathcal{T}^0(t_k, x)) = y\}$ by (up to a multiplicative constant) the probability $\mathbb{P}\{\Pi_\infty(\mathcal{T}(t_k, x)) = y\}$. We set

$$\Sigma(t_{k+1}, x) = \sigma(x, \bar{u}(t_{k+1}, x)), \qquad \mu(t_{k+1}, x) = b(x, \bar{u}(t_{k+1}, x), \hat{v}(t_k, x)).$$

Put $\|\Sigma(\text{resp. } \mu)\|_\infty = \sup_{k \in \{0,\dots,N\}}[\sup_{x \in \mathcal{C}_k} |\Sigma(\text{resp. } \mu)(t_k, x)|]$. From Assumption (A) and Propositions 6.1 and 6.2 (boundedness of $\bar{u}$ and $h^{1/2}\hat{v}$), $\|\Sigma\|_\infty + h^{1/2}\|\mu\|_\infty \le C$.

*Step two*: *Proof of Theorem* 7.6, *comparison of kernels.* The proof of the following proposition relies on a standard Gaussian change of variable and rather tedious computations (the detailed proof is given in Section 7.3 in the electronic version [9]):

PROPOSITION 7.7. *There exists a constant* $C_{7.7} > 0$ *such that, for every* $y \in \mathcal{C}_\infty$,

$$\mathbb{P}\{\Pi_\infty(\mathcal{T}^0(t_k, x)) = y\} \leq \alpha_k(y) + \beta(y)(\eta_k + \mathbb{P}^{1/2}\{\Pi_\infty(\mathcal{T}(t_k, x)) = y\}),$$

*where*

$$\alpha_k(y) \equiv \mathbb{P}\{|\Sigma(t_{k+1}, x)g(\Delta B^k) - y|_\infty \leq \delta/2,$$
$$|g(\Delta B^k) - \Delta B^k|_\infty > \delta/(2\|\Sigma\|_\infty)\},$$
$$\beta(y) \equiv C_{7.7}\delta^{d/2}h^{-d/4}\exp[-C_{7.7}^{-1}h^{-1}|y|^2],$$
$$\eta_k \equiv \mathbb{P}^{1/2}\{|g(\Delta B^k) - \Delta B^k|_\infty > \delta/(4\|\Sigma\|_\infty)\}.$$

*In the above expression, for all* $z \in \mathbb{R}^d$, $|z|_\infty \equiv \max_{i \in \{1,\dots,d\}} |z_i|$.

From Proposition 6.2, $h^{1/2}\bar{v}$ is bounded by a known constant. Denote by RHS$(X_{t_k}, 7.13)$ the r.h.s. in (7.13) and by $\Gamma(h, C)$ the sum $\sum_{y \in \mathcal{C}_\infty} \exp[-C^{-1} \times h^{-1}|y|^2]$. Owing to Proposition 7.7 and (7.12), we then get

$$\sum_{y \in \mathcal{C}_\infty} [|(\bar{v} - v)(t_{k+1}, \Pi_{k+1}(x + y))|\mathbb{P}\{\Pi_\infty(\mathcal{T}^0(t_k, x)) = y\}]$$

(7.14)
$$\leq Ch^{-1/2}\mathbb{P}\{|g(\Delta B^k) - \Delta B^k|_\infty > \delta/(2\|\Sigma\|_\infty)\}$$
$$+ C\delta^{d/2}h^{-d/4-1/2}\mathbb{P}^{1/2}\{|g(\Delta B^k) - \Delta B^k|_\infty > \delta/(4\|\Sigma\|_\infty)\}$$
$$\times \Gamma(h, C)$$
$$+ C[\delta^d h^{-d/2}\Gamma(h, C)]^{1/2}[\text{RHS}(x, 7.13)]^{1/2}$$
$$\equiv T(1) + T(2) + T(3).$$

Due to (3.8) and to the Bienaymé–Chebyshev inequality, $T(1) \leq Ch^{p/2-1/2}$ $\delta^{-p}M^{-p/d}$. Thanks again to (3.8) (applied to the exponent $2p$), $T(2) \leq Ch^{p/2-d/4-1/2}\delta^{-p+d/2}M^{-p/d}\Gamma(h, C) = C\mathcal{E}(\text{gradient}, p)(\delta h^{-1/2})^d\Gamma(h, C)$.

Note now from (7.13) that

$$T(3) = C[\delta^d h^{-d/2}\Gamma(h, C/2)]^{1/2}\mathbb{E}[|(\bar{v} - v)(t_{k+1}, \Pi_{k+1}(x + \mathcal{T}(t_k, x)))|^2]^{1/2}.$$

A standard comparison with a Gaussian integral yields $(\delta h^{-1/2})^d\Gamma(h, C) \leq C''$. Plugging the different estimates of $T(1)$, $T(2)$ and $T(3)$ in (7.14), we complete the proof of Theorem 7.6 [recall again that $h^{-1}\delta^2$ is small to dominate $T(1)$ by $\mathcal{E}(\text{gradient}, p)$].

**8. Proof. Third step: Gronwall's lemma.** Here is the final step of the proof of Theorems 4.1, 4.2 and 4.3.

8.1. *Proof of Theorem* 4.1, *infinite grids.* We first explain how to derive Theorem 4.1 from Theorem 7.3 when $\rho = +\infty$, that is, when $\tau_\infty = +\infty$ a.s. In this framework, the term $\mathcal{E}^2(\text{trunc})$ in $\mathcal{E}^2(\text{global})$ reduces to 0. The general case is detailed in the next subsection. For infinite grids, for $\eta$ and $h$ small enough, we obtain from Theorem 7.3 and from the equality $\bar{v}(T,x) = v(T,x)$, for all $x \in \mathcal{C}_N$,

$$(8.1) \quad |(\bar{u} - u)(0, x_0)|^2 \le C \left[ \mathcal{E}^2(\text{global}) + h \sum_{j=0}^{N} \sup_{x \in \mathcal{C}_j} |(\bar{u} - u)(t_j, x)|^2 \right].$$

As usual in BSDE theory, the estimate (8.1) holds actually for any starting point $(t_k, x)$, $0 \le k \le N$, $x \in \mathcal{C}_k$. Hence, there is no difficulty to apply Gronwall's lemma [at least for $h$ small, as in (6.7)] and to complete the proof of Theorem 4.1 when $\rho = +\infty$.

8.2. *Proof of Theorem* 4.1, *general case.* We now turn to the case of truncated grids. Generally speaking, most of the approach given in the former subsection still applies in the general framework. It is, however, impossible to mimic word for word the arguments given above and we need to refine the previous Gronwall argument.

*First step.* We first aim to get rid of the difference $\bar{v} - v$ appearing in the new r.h.s. in Theorem 7.3. Due to the functions $(\mathbf{1}_{\{t_{j-1} < \tau_\infty\}})_{j=1,\dots,N}$, the machinery used in the previous subsection does not apply. To overcome this difficulty, we write $\{t_{j-1} < \tau_\infty\} = \{t_j < \tau_\infty\} \cup \{t_j = \tau_\infty\}$. Indeed, since $\bar{v}(T,x) = v(T,x)$ for $x \in \mathcal{C}_N$ and $h|\bar{v} - v|^2$ is bounded (see Theorem 2.1 and Proposition 6.2), we obtain for $\eta$ and $h$ small enough

$$|(\bar{u} - u)(0, x_0)|^2 + C^{-1} h \sum_{j=1}^{N} \mathbb{E}[|(\bar{v} - v)(t_{j-1}, X_{t_{j-1}})|^2 \mathbf{1}_{\{t_{j-1} < \tau_\infty\}}]$$

$$\le C \left[ \mathbb{P}\{\tau_\infty < +\infty\} + \mathcal{E}^2(\text{global}) \right.$$

$$(8.2) \qquad + h \sum_{j=1}^{N} \mathbb{E}[|(\bar{u} - u)(t_j, X_{t_{j-1}})|^2 \mathbf{1}_{\{t_{j-1} < \tau_\infty\}}]$$

$$\left. + h \sum_{j=2}^{N} \mathbb{E}[|(\bar{u} - u)(t_{j-1}, X_{t_{j-1}})|^2 \mathbf{1}_{\{t_{j-1} < \tau_\infty\}}] \right].$$

Even though we employed $\mathcal{E}^2(\text{global})$ for notational convenience, we mention carefully that the origin of the term $\mathcal{E}^2(\text{trunc})$ has not been explained yet. It is in the following lines.

*Second step.* Note that (8.2) still holds if $X$ starts at a given time $t_i$, $i \in \{0, \ldots, N\}$, from an $\mathcal{F}_{t_i}$-measurable and square integrable random variable $\xi$ with values in $\mathcal{C}_i$. In such a case, (8.2) still holds with $(0, x_0)$ replaced by $(t_i, \xi)$, $X_{t_j}$ by $X_{t_j}^{t_i, \xi}$ and $\tau_\infty$ by $\tau_\infty^{t_i, \xi}$ [the superscript $(t_i, \xi)$ denotes the initial condition of the process $X$]. Due to the shift between $t_{j-1}$ and $t_j$ in the r.h.s., there is no possible choice of $\xi$ to recover the same form of terms in the left and right-hand sides. In particular, Gronwall's lemma does not apply at this stage of the proof. Note, in fact, that the same problem occurred in Section 8.1: this was the reason why the supremum was taken in the r.h.s. of (8.1).

In the current frame, taking the supremum over $x \in \mathcal{C}_i$ in (8.2) induces a new term, namely, $\sup_{x \in \mathcal{C}_i} \mathbb{P}\{\tau_\infty^{t_i, x} < +\infty\}$. Unfortunately, for $x$ close to the boundary of the grid $\mathcal{C}_i$, the underlying probability is far from being small. In particular, there is no hope to prove Theorem 4.1 in the case $\rho < +\infty$ with the arguments used in Section 8.1.

*Strategy.* Our strategy then consists in applying (8.2) to a suitable choice of $\xi$. We then have to estimate the probability $\mathbb{P}\{\tau_\infty^{t_i, \xi} < +\infty\}$ for a random initial condition $(t_i, \xi)$, $\xi \in L^2(\Omega, \mathcal{F}_{t_i}, \mathbb{P})$ with values in $\mathcal{C}_i$. To this end, we need to control efficiently the tails of the variables $(X_{t_j \wedge \tau_\infty^{t_i, \xi}}^{t_i, \xi})_{i \leq j \leq N}$. Since the drift $b$ is not bounded, a natural approach consists in estimating the $L^2$ norms of these variables.

LEMMA 8.1 ($L^2$ *control of the process* $X$). *For all* $k \in \{0, \ldots, N\}$, *put* $\tau_k = \tau_\infty \wedge t_k$. *Then, there exists a constant* $C_{8.1}$ *such that, for all* $i \in \{0, \ldots, N\}$ *and* $\xi \in L^2(\Omega, \mathcal{F}_{t_i}, \mathbb{P})$ *with values in* $\mathcal{C}_i$,

$$\forall k \in \{i, \ldots, N\}, \qquad \mathbb{E}[|X_{\tau_k}^{t_i, \xi}|^2] \leq C_{8.1}[\mathbb{E}|\xi|^2 + 1 + \mathcal{E}^2(\text{space}) + \mathcal{E}^2(\text{gradient}, p)].$$

PROOF (SKETCH). We remove the superscript $(t_i, \xi)$ in the writing of $X$. Then

$$X_{\tau_k} = \xi + \sum_{j=i}^{k-1} [\mathcal{T}(t_j, X_{t_j}) \mathbf{1}_{\{t_j < \tau_\infty\}}]$$

$$+ \sum_{j=i}^{k-1} [(\Pi_{j+1}(X_{t_j} + \mathcal{T}(t_j, X_{t_j})) - X_{t_j} - \mathcal{T}(t_j, X_{t_j})) \mathbf{1}_{\{t_{j+1} < \tau_\infty\}}]$$

(8.3)

$$+ \sum_{j=i}^{k-1} [(\Pi_{j+1}(X_{t_j} + \mathcal{T}(t_j, X_{t_j})) - X_{t_j} - \mathcal{T}(t_j, X_{t_j}))\mathbf{1}_{\{t_{j+1}=\tau_\infty\}}]$$

$$\equiv \xi + S(1) + S(2) + S(3).$$

The term $S(2)$ corresponds to a standard projection error. Thus, $\mathbb{E}[|S(2)|^2] \leq \delta^2(k-i)^2 \leq C\mathcal{E}^2(\text{space})$. For $S(3)$, Lemma 6.4 and Young's inequality yield

$$
\begin{aligned}
(8.4) \quad & \mathbb{E}[|S(3)|^2] \leq C\delta^2 + C\sum_{j=i}^{k-1} \mathbb{E}[|\mathcal{T}(t_j, X_{t_j})|^2 \mathbf{1}_{\{t_{j+1}=\tau_\infty\}}] \\
& \leq Ch^2\mathcal{E}^2(\text{space}) + C\mathbb{P}\{\tau_\infty < +\infty\} + C\sum_{j=i}^{k-1} \mathbb{E}[|\mathcal{T}(t_j, X_{t_j})|^4].
\end{aligned}
$$

From Propositions 6.1 and 6.2, we can prove that $\mathbb{E}[|\mathcal{T}(t_j, X_{t_j})|^4] \leq Ch^2$. We finally deduce ($h$ being small) $\mathbb{E}[|S(3)|^2] \leq C[\mathbb{P}\{\tau_\infty < +\infty\} + \mathcal{E}^2(\text{time}) + \mathcal{E}^2(\text{space})]$.

Deal now with $S(1)$. Thanks to Propositions 6.1 and 6.2, we estimate the drift, and thanks to the independence of the Brownian increments, we bound the martingale part. From Assumption (A), there exists a constant $C$ such that

$$\mathbb{E}[|S(1)|^2] \leq Ch(k-i)\left[1 + h\sum_{j=i}^{k-1} \mathbb{E}[|\hat{v}(t_j, X_{t_j})|^2 \mathbf{1}_{\{t_j<\tau_\infty\}}]\right].$$

Apply now Propositions 7.2 and 6.2, and derive that $\mathbb{E}[|S(1)|^2] \leq C[1 + \mathcal{E}^2(\text{gradient}, p)]$. $\square$

*Estimate of the probability of hitting the boundary.* Thanks to the previous lemma, we are now able to estimate the probability $\mathbb{P}\{\tau_\infty^{t_i,\xi} < +\infty\}$, with $(i,\xi)$ as in Lemma 8.1. Indeed, $\{\tau_\infty^{t_i,\xi} < +\infty\} \subset \{|X_{\tau_N}^{t_i,\xi}|_\infty + \delta \geq R + \rho\}$. Thanks to the Bienaymé–Chebyshev inequality and to Lemma 8.1 (with $k = N$), we get

$$
\begin{aligned}
(8.5) \quad & \mathbb{P}\{\tau_\infty^{t_i,\xi} < +\infty\} \\
& \leq C[(R+\rho)^{-2}\mathbb{E}[|\xi|^2] + \mathcal{E}^2(\text{space}) + \mathcal{E}^2(\text{trunc}) + \mathcal{E}^2(\text{gradient}, p)].
\end{aligned}
$$

Plug now (8.5) into (8.2) to obtain

$$
\begin{aligned}
(8.6) \quad & \mathbb{E}[|(\bar{u} - u)(t_i, \xi)|^2] \\
& \leq C\left[(R+\rho)^{-2}\mathbb{E}[|\xi|^2] + \mathcal{E}^2(\text{global})\right.
\end{aligned}
$$

$$+ h \sum_{j=i+1}^{N} \mathbb{E}[|(\bar{u}-u)(t_j, X_{t_{j-1}}^{t_i, \xi})|^2 \mathbf{1}_{\{t_{j-1} < \tau_\infty^{t_i, \xi}\}}]$$

$$+ h \sum_{j=i+2}^{N} \mathbb{E}[|(\bar{u}-u)(t_{j-1}, X_{t_{j-1}}^{t_i, \xi})|^2 \mathbf{1}_{\{t_{j-1} < \tau_\infty^{t_i, \xi}\}}]\Bigg].$$

*A refined Gronwall argument.* The key idea is to find by induction a sequence of constants $c_i(1)$, $c_i(2)$, $i \in \{0, \ldots, N\}$, such that, for any $\xi \in L^2(\Omega, \mathcal{F}_{t_i}, \mathbb{P})$ with values in $\mathcal{C}_i$,

(8.7)
$$\mathbb{E}[|(\bar{u}-u)(t_i, \xi)|^2]$$
$$\leq c_i(1)\mathcal{E}^2(\text{global}) + c_i(2)(R+\rho)^{-2}\mathbb{E}[|\xi|^2].$$

Thanks to Lemma 8.1, we are able to build two sequences $c_i(1)$ and $c_i(2)$, $i \in \{0, \ldots, N\}$, satisfying (8.7) and uniformly bounded by a constant $C$. Choosing $i = 0$ and $\xi = x_0 \in \mathcal{C}_0$, we then complete the proof of Theorem 4.1. The explicit construction of $c_i(1)$ and $c_i(2)$, $i \in \{0, \ldots, N\}$, is given in the electronic version [see (8.12) in there].

8.3. *Proofs of Theorems 4.2 and 4.3.* We turn to the proof of Theorems 4.2 and 4.3. The initial condition of the process $X$ is given by $X_0 = x_0$, $x_0 \in \mathcal{C}_0$, as in (3.12).

PROOF OF THEOREM 4.2. From inequalities (8.2) (deriving from the stability theorem), (8.5) (probability of hitting the boundary of the grids) and (8.7) [estimate of $\bar{u} - u$, recall that $c_j(1)$, $c_j(2)$, $j \in \{0, \ldots, N\}$, are uniformly bounded], Theorem 4.2 holds with $v(t_i, X_{t_i})\mathbf{1}_{\{t_i < \tau_\infty\}}$ instead of $v(t_i, X_{t_i})$. Since $v$ is bounded (see Theorem 2.1) and since the probability of hitting the boundaries of the grids is controlled [see again (8.5)], we easily complete the proof. □

PROOF OF THEOREM 4.3. It just remains to study the convergence of $(X_{t_k}, Y_{t_k}, Z_{t_k})_{0 \leq t_k \leq \tau_\infty \wedge T}$ toward the solution $(U, V, W)$ of (E). Thanks to the Lipschitz properties of $b$ and $\sigma$, we first deduce by standard computations (see, e.g., the proof of Lemma 8.1) the analogue of Proposition 7.1.

PROPOSITION 8.2. *There exists a constant $C_{8.2}$ s.t., for $k \in \{1, \ldots, N\}$,*

$$\mathbb{E}|X_{\tau_k} - U_{\tau_k}|^2$$
$$\leq C_{8.2}\Bigg[\mathbb{P}\{\tau_\infty < +\infty\} + \mathcal{E}^2(\text{global})$$

(8.8)

$$+ h \sum_{j=0}^{k-1} [\mathbb{E}[\mathbf{1}_{\{t_j < \tau_\infty\}}(|X_{t_j} - U_{t_j}|^2 + |(\bar{u} - u)(t_{j+1}, X_{t_j})|^2$$

$$+ |(\hat{v} - v)(t_j, X_{t_j})|^2)]].$$

Recall now from Proposition 7.2 (estimate of $\hat{v} - v$), Theorem 4.2 ($L^2$ estimate of $\bar{v} - v$) and (8.5) (probability of hitting the boundary of the grids):

(8.9)
$$h \sum_{j=0}^{k-1} \mathbb{E}[|(\hat{v} - v)(t_j, X_{t_j})|^2 \mathbf{1}_{\{t_j < \tau_\infty\}}]$$

$$\leq C \Bigg[ \mathcal{E}^2(\text{time}) + \mathcal{E}^2(\text{space}) + \mathcal{E}^2(\text{gradient}, p)$$

$$+ h \sum_{j=1}^{k} \mathbb{E}[|(\bar{v} - v)(t_j, X_{t_j})|^2 (\mathbf{1}_{\{t_j < \tau_\infty\}} + \mathbf{1}_{\{t_j = \tau_\infty\}})] \Bigg]$$

$$\leq C[\mathcal{E}^2(\text{global}) + \mathbb{P}\{\tau_\infty < +\infty\}] \leq C\mathcal{E}^2(\text{global}).$$

Apply now inequality (8.7) (estimate of $\bar{u} - u$) and (8.9) to (8.8) and deduce from Gronwall's lemma that $\sup_{k \in \{0,\dots,N\}} \mathbb{E}|X_{\tau_k} - U_{\tau_k}|^2 \leq C\mathcal{E}^2(\text{global})$. Finally, according to Theorem 2.1, to Theorem 4.2 ($L^2$ estimate of $\bar{v} - v$) and to (8.7), we deduce the following intermediate estimate:

(8.10)
$$\sup_{k \in \{0,\dots,N\}} \mathbb{E}[|X_{\tau_k} - U_{\tau_k}|^2 + |Y_{\tau_k} - V_{\tau_k}|^2]$$

$$+ h \sum_{j=0}^{N-1} \mathbb{E}[|Z_{t_j} - W_{t_j}|^2 \mathbf{1}_{\{t_j < \tau_\infty\}}] \leq C\mathcal{E}^2(\text{global}).$$

Applying Doob's inequality, we derive the same bound but with the *supremum* inside the expectation. It finally remains to prove the same result, but with $(U_{t_k}, V_{t_k}, W_{t_k})_{0 \leq k \leq N}$ instead of $(U_{\tau_k}, V_{\tau_k}, W_{t_k} \mathbf{1}_{\{t_k < \tau_\infty\}})_{0 \leq k \leq N}$. Since the same arguments apply for $V$ and $W$, we just detail the case of $U$. Note indeed that, for every $k \in \{0, \dots, N\}$,

$$\sup_{k \in \{0,\dots,N\}} |X_{\tau_k} - U_{t_k}|^2 \leq C \sup_{k \in \{0,\dots,N\}} |X_{\tau_k} - U_{\tau_k}|^2 + C \sup_{k \in \{0,\dots,N\}} |U_{\tau_k} - U_{t_k}|^2.$$

Thanks to the Burkholder, Davis and Gundy inequalities, it is readily seen that

$$\mathbb{E}\Big[\sup_{k \in \{0,\dots,N\}} |U_{\tau_k} - U_{t_k}|^2\Big] \leq C\mathbb{E}[(t_N - \tau_\infty)\mathbf{1}_{\{\tau_\infty < +\infty\}}] \leq CT\mathbb{P}\{\tau_\infty < +\infty\}.$$

Referring to (8.5), we easily complete the proof of Theorem 4.3.  □

**9. Conclusion.** As a conclusion, we first give in Section 9.1 further comments on Theorem 4.1 and compare, in particular, the global error with the one obtained by Douglas, Ma and Protter [10]. We then give some easy extensions in Section 9.2. Finally, we detail in Section 9.3 the technical difficulties associated with the natural algorithm (3.9)–(3.10).

9.1. *Comments and comparisons with other methods.* We discuss in this subsection the total complexity and the rate of convergence of Algorithm 3.1.

*Complexity of the algorithm.* Note first that the order of the total complexity of the algorithm is $h^{-1} \times M \times (2\delta^{-1}(\rho + R))^d$.

*Rate of convergence.* Recall also that the global error of the algorithm is given by Theorem 4.1. Comparing with the results in [10], this global error is worse in our case. There are two reasons to explain this difference. The first one does not depend on the algorithm, but is a consequence of our working assumptions. Indeed, under suitable smoothness properties of the coefficients $b, f, \sigma$ and of the solution $u$, standard Itô developments in $D(1)$ (see Lemma 7.4) would lead to $\mathcal{E}^2(\text{time}) = h^2$ as in [10].

At the opposite, the second reason for which the global error is worse, in our case, depends on the specific structure of the algorithm. Indeed, our choice to avoid linear interpolation procedures induces a rather large projection error $\mathcal{E}^2(\text{space})$. To reach a term of order one with respect to $h$ for $\mathcal{E}^2(\text{space})$, we then need to take $\delta \equiv h^{3/2}$. This choice is far from being satisfactory and highly increases the complexity when the dimension $d$ increases. Intuitively, there is no specific reason for such a relationship between $\delta$ and $h$: as explained in Section 4.1, $\delta$ has just to be small in front of $h$ to take into account the influence of the drift $b$ at the local level. For this reason, we aim to study in further investigations the convergence analysis of the algorithm when using a suitable "smooth" interpolation operator instead of a rough projection mapping. This point is discussed in a detailed way in the next subsection.

*Further comments on errors.* To conclude this subsection, we investigate the three last error terms, $\mathcal{E}(\text{trunc})$, $\mathcal{E}(\text{quantiz})$ and $\mathcal{E}(\text{gradient}, p)$.

The truncation error decays linearly when the grid size increases. This control may seem rather poor to the reader. Recall indeed that $\mathcal{E}(\text{trunc})$ appears, up to the discretization procedure, as the probability that a diffusion process leaves a given bounded set. In the case of elliptic diffusions with bounded coefficients, it is well known that this probability decays exponentially fast as the size of the underlying set increases. Recall in this frame from Theorem 2.1 that the coefficients of the elliptic diffusion $U$ are bounded. Note, however, that this rough argument fails in the discretized setting since there is no a priori sharp estimate of the approximate gradient $\bar{v}$ and thus of the associated approximate drift. This explains why our strategy to estimate $\mathcal{E}(\mathrm{trunc})$ lies on the Bienaymé–Chebyshev inequality and, thus, provides the current form given by Theorem 4.1. Similar techniques could yield a polynomial decay for every $q \geq 1$, the constant of the theorem being an increasing function of $q$ [see Lemma 8.1 and (8.5)].

Note finally that the errors associated to the quantization procedure, $\mathcal{E}^2(\mathrm{quantiz})$, and to the probabilistic approximation of the gradient, $\mathcal{E}^2(\mathrm{gradient}, p)$, are explicitly controlled in terms of $M$, $h$ and $\delta$. They emphasize the price to pay to weaken the assumptions: we have to assume that the quantization grid is rather small compared to the spatial discretization one. Obviously, this increases the number of elementary operations of the algorithm and, thus, its total complexity. However, this does not affect so much the discretization procedure of the Gaussian law itself since quantization grids can be computed once for all apart from the implementation procedure of the algorithm.

### 9.2. *Extensions and further investigations.* We now discuss some possible extensions of our work.

*Interpolation procedure.* As stated later in this subsection, we first investigate the assets and liabilities of a smooth interpolation procedure. One of the main advantages of the spatial discretization proposed in Section 3.4, and then used in Algorithm 3.1, lies in its simplicity of implementation. However, from a purely mathematical point of view, this procedure may be rather awkward since it ignores more or less the deep smoothness of the true solution $u$.

Note in this framework that the function $\Pi_\infty$ may be seen as an operator acting on functions from $\mathbb{R}^d$ into $\mathbb{R}$. For such a function, the operator provides a rough interpolation of order 0 depending on the values of the function on the spatial mesh $\mathcal{C}_\infty$. As mentioned above, this interpolation procedure does not preserve the smoothness properties of the underlying function: in any cases, except if the function is constant, the interpolation procedure induces jumps of size of order $\delta$. As a consequence, the distance between the function and the interpolated one is also of order $\delta$.

A relevant strategy would consist in replacing the projection $\Pi_\infty$ by a smoother interpolation operator. In our framework, an interpolation operator is said to be smooth if the distance between a given function $\ell$ and the interpolated one decreases with the regularity order of $\ell$. For example, in dimension 1, the linear interpolation operator,

$$\ell \mapsto (x \mapsto \delta^{-1}(\delta + \delta \lfloor \delta^{-1}x \rfloor - x)\ell(\delta \lfloor \delta^{-1}x \rfloor) + \delta^{-1}(x - \delta \lfloor \delta^{-1}x \rfloor)\ell(\delta \lfloor \delta^{-1}x \rfloor + \delta)),$$

maps a $\mathcal{C}^2(\mathbb{R}, \mathbb{R})$ function into a piecewise smooth function and the distance between them is of order $\delta^2$.

Algorithm 3.1 can be written with respect to this new choice, but we also believe that the proof would be more difficult to detail. Moreover, smooth interpolation procedures in higher dimension slow down the running of the underlying algorithm.

*Weakening assumption.* Note to conclude this subsection that some assumptions could be weakened. First, Theorem 2.1 still holds if $b$ and $f$ are just Hölder in $x$: in such a case, usual estimates of the gradient of $u$ hold and Schauder's theory still applies. In particular, the reader can verify that Theorems 4.1 and 4.2 are still valid in this case (but Theorem 4.3 given in Section 4.3 is not).

Moreover, Algorithm 3.1 still converges if $b, f$ and $\sigma$ depend on $t$ in a Hölder way.

Finally, the following extension is conceivable. For $H \in C^{1+\alpha}$, $\alpha \in ]0, 1[$, the partial derivatives of order two of $u$ have an integrable singularity in the neighborhood of $T$. In this frame, it would be interesting to adapt the Gronwall arguments given in Section 8.

9.3. *Justification of Algorithm* 3.1. We finally explain why we are not able to show the convergence of Algorithm (3.9)–(3.10).

*Convergence of algorithm* (3.9)–(3.10). Recall that the main difference between the algorithm (3.9)–(3.10) and Algorithm 3.1 lies in the definition of the forward transitions. Indeed, in the algorithm (3.9)–(3.10),

$$\mathcal{T}(t_k, x) \equiv b(x, \bar{u}(t_{k+1}, x), \bar{v}(t_{k+1}, x))h + \sigma(x, \bar{u}(t_{k+1}, x))g(\Delta B^k),$$

$$X_0 \equiv x_0, \qquad \forall k \in \{0, \ldots, N-1\}, \qquad X_{t_{k+1}} = \Pi_{k+1}(X_{t_k} + \mathcal{T}(t_k, t_k)).$$

Unfortunately, in this case, the well-known BSDE machinery fails under Assumption (A). At first sight, this could seem rather amazing. Indeed, recall that very strong a priori estimates of the solution $u$ and of its partial derivatives hold in our framework. In particular, we could expect the discretization procedure of $u$ and of its gradient to converge under such smoothness properties.

The main difficulty encountered to establish the convergence of the algorithm (3.9)–(3.10) appears in Section 7. The lack of a priori controls of the regularity of $\bar{u}$ and $\bar{v}$ makes the stability strategy fruitless. Note, indeed, that inequality (7.1) becomes in the frame of the indicated algorithm

$$|(\bar{u} - u)(0, x_0)|^2 + C_{7.1}^{-1} h \sum_{j=1}^{N} \mathbb{E}[|(\bar{v} - v)(t_{j-1}, X_{t_{j-1}})|^2 \mathbf{1}_{\{t_{j-1} < \tau_\infty\}}]$$

$$\leq C_{7.1} \Bigg[ \mathbb{P}\{\tau_\infty < +\infty\} + \mathcal{E}^2(\text{time}) + \mathcal{E}^2(\text{space}) + \mathcal{E}^2(\text{quantiz})$$

(9.1)
$$+ \eta^{-1} h \sum_{j=1}^N \mathbb{E}[|(\bar{u} - u)(t_j, X_{t_{j-1}})|^2 \mathbf{1}_{\{t_{j-1} < \tau_\infty\}}]$$

$$+ \eta^{-1} h \sum_{j=1}^N \mathbb{E}[|(\bar{u} - u)(t_{j-1}, X_{t_{j-1}})|^2 \mathbf{1}_{\{t_{j-1} < \tau_\infty\}}]$$

$$+ (\eta + h) h \sum_{j=1}^N \mathbb{E}[|(\bar{v} - v)(t_j, X_{t_{j-1}})|^2 \mathbf{1}_{\{t_{j-1} < \tau_\infty\}}] \Bigg].$$

Inequalities (7.1) and (9.1) just differ in the last term: $\hat{v}(t_{j-1}, X_{t_{j-1}})$ becomes $\bar{v}(t_j, X_{t_{j-1}})$. Note that to be complete a similar shift occurs in $v$ but, due to Theorem 2.1, it can be removed without any difficulties. To apply the strategy used in Section 7, and, in particular, to derive an equivalent of Theorem 7.3 from (9.1), we then need to investigate the regularity in space of $\bar{v}$. According to the definition of $\bar{v}$, a first step then consists in studying the regularity in space of $\bar{u}$.

*Lipschitz control of $\bar{u}$.* Note that the natural strategy to control the oscillations of $\bar{u}$ would consist in applying the usual FBSDE machinery to the triples $(X^{t_k,x}, Y^{t_k,x}, Z^{t_k,x})$ and $(X^{t_k,y}, Y^{t_k,y}, Z^{t_k,y})$ for $k \in \{0, \ldots, N-1\}$ and $x, y \in \mathcal{C}_k$. Of course, superscripts $(t_k, x)$ and $(t_k, y)$ denote the initial conditions of the Markov process $X$.

Nevertheless, we are not able to apply the strategies used in [6, 7] to derive from the forward–backward writing local and global estimates of the discrete gradient of $\bar{u}$. There are two reasons to explain this failure.

First, the rough projection mapping chosen induces an irreducible error greater than $\delta$ when estimating the difference between $\bar{u}(t_k, x)$ and $\bar{u}(t_k, y)$ in function of the parameters deriving from Assumption (A). The strategy to overcome this difficulty is well known: the projection mapping has to be replaced by a smoother interpolation operator.

Second, any probabilistic strategy to estimate the Lipschitz constant of $\bar{u}$ in $x$ such as the one exposed in [6] leads one way or another to the same difficulty as the one encountered to apply the stability procedure to the algorithm (3.9)–(3.10). More precisely, studying the difference between the triples $(X^{t_k,x}, Y^{t_k,x}, Z^{t_k,x})$ and $(X^{t_k,y}, Y^{t_k,y}, Z^{t_k,y})$, for $k \in \{0, \ldots, N-1\}$ and $x, y \in \mathcal{C}_k$, leads to investigate the regularity of $\bar{v}$. In short, one needs to estimate first the regularity of $\bar{v}$ to derive the one of $\bar{u}$. Intuitively, it is well understood that this is hopeless.

**Acknowledgment.** We would like to thank Professor Vlad Bally for his thoughtful attentive comments.

UNIVERSITÉ PARIS VII
UFR DE MATHÉMATIQUES
CASE 7012
2 PLACE JUSSIEU
75251 PARIS CEDEX 05
FRANCE
E-MAIL: delarue@math.jussieu.fr
         menozzi@ccr.jussieu.fr